\documentclass{article}

\usepackage{arxiv}

\usepackage[utf8]{inputenc} 
\usepackage[T1]{fontenc}    
\usepackage{hyperref}       
\usepackage{url}            
\usepackage{booktabs}       
\usepackage{amsfonts}       
\usepackage{nicefrac}       
\usepackage{microtype}      
\usepackage{lipsum}		
\usepackage{graphicx}
\usepackage{natbib}
\usepackage{doi}
\usepackage{amsmath}
\usepackage{algorithm}
\usepackage{algpseudocode}
\usepackage{amssymb}

\title{Contribution of expert aggregation to temperature prediction Part II: second-order bounds with sleeping experts.}

\author{ {Leo ~Pfitzner} \\
	Meteo France\\
	\texttt{leo.pfitzner@meteo.fr} \\
	\And
	{Olivier~Wintenberger} \\
	Sorbonne University\\
	\texttt{olivier.wintenberger@sorbonne-universite.fr} \\
    \And
	{Olivier~Mestre} \\
	Meteo France\\
	\texttt{olivier.mestre@meteo.fr}
}

\renewcommand{\shorttitle}{}

\hypersetup{
pdftitle={A template for the arxiv style},
pdfsubject={q-bio.NC, q-bio.QM},
pdfauthor={David S.~Hippocampus, Elias D.~Striatum},
pdfkeywords={First keyword, Second keyword, More},
}

\begin{document}
\maketitle

\begin{abstract}
In this paper, we improve on the temperature predictions made with online Expert Aggregation (EA) \citep{cesa-bianchi_prediction_2006} in Part I \citep{}. In particular, we make the aggregation more reactive, whilst improving the root mean squared error and reducing the number of large errors. We have achieved this by using the Sleeping Expert Framework (SEF) \citep{freund_using_1997,devaine_forecasting_2013}, which allows a more efficient use of biased experts (bad on average but which may be good at some point). To deal with the fact that, unlike in \citet{devaine_forecasting_2013}, we do not know in advance when to use these biased experts, we resorted to Gradient Boosted Regression Trees (GBRT) \citep{chen_xgboost_2016} and provide regret bounds against sequences of experts \citep{mourtada_efficient_2017} which take into account this uncertainty. We applied the SEF with GBRT in a fully online way on Bernstein Online Aggregation (BOA) \citep{wintenberger_stochastic_2024}, an adaptive aggregation with second-order regret bounds, which had the best results in Part I.
\end{abstract}


\section{Introduction}
Temperature forecasts are crucial in many sectors, such as agriculture, road management and tourism. Today, numerous forecasts are available, making it desirable to integrate them for a more accurate prediction.

This can be achieved online using Expert Aggregation (EA), with theoretical guarantees, as presented in Part I. We recall that in the EA framework the different predictions/models to be aggregated are also called "experts". An overview of EA can be found in \citet{cesa-bianchi_prediction_2006}.

In Part I we introduced and compared some EAs and explained how they can improve the average temperature predictions of the raw and Post-Processed Models (PPM) of the Numerical Weather Prediction (NWP) models by competing with the best expert - or the best convex combination of experts - in hindsight. Avoiding large (consecutive) forecast errors may however be even more important than having good forecasts on average.

The main problem identified in Part I is that the EA rapidly loses the ability to use biased experts effectively because they are bad on average. Yet, these biased experts could help to avoid large errors because extreme events are often only well predicted by the low or high quantiles of an ensemble forecast \citep{taillardat_calibrated_2016,taillardat_research_2020}.

The addition of a sliding window in Part I made the EA more reactive and more able to use any expert at any time. However, the drawback of a sliding window is that it reduces performance on average. 

Hence, it remains to make the EA more reactive and able to use the biased experts efficiently without adding noise in a changing environment where the best expert can change over time. As in Part I, we recall that we are trying to improve predictions that are already good on average. Here, however, the task is even more challenging, since we aim to further improve the predictions obtained in Part I of the EA itself, which already largely outperformed the uniform aggregation of the models and the best expert in hindsight including PPMs. Therefore, in this paper, unlike in Part I, we focus on aggregation strategies which are designed to handle changing environments. Thus, we address a question left open by Part I: how can an otherwise strong online EA framework make
better use of experts that are poor on average but useful in short-lived or difficult situations? 

There are three main ways to deal with changing environments in EA. The first way is to use so-called "meta-aggregations" which aim to minimize regret over any possible time interval \citep{hazan_adaptive_2007,jun_improved_2017} or, alternatively, to learn the transition probabilities of the changing environment \citep{monteleoni_tracking_2011}.

The second way, introduced by \citet{herbster_tracking_1998} with its Fixed Share (FS) strategy (and generalized by \citet{mourtada_efficient_2017}), deals with changing environments by competing with sequences of experts. This framework has been used for electricity consumption by \citet{devaine_forecasting_2013} and \citet{vyugin_online_2019}.

The third way - which we have used in this paper - is the Specialized (or Sleeping) Experts Framework (SEF), introduced by \citet{freund_using_1997} and \citet{blum_empirical_1997} and then used for electricity consumption forecasting in \citet{devaine_forecasting_2013}. This framework tries to take advantage of the fact that so-called specialized/sleeping experts can be good at certain events and bad at others.

In \citet{devaine_forecasting_2013} and \citet{vyugin_online_2019}, the authors used sleeping experts which were specialized in predicting electricity consumption for certain periods (e.g. holidays, winter, summer). So it is known in advance when these experts should make good predictions.

To the best of our knowledge, the SEF framework has never been used with second-order EAs as it is in this article, nor in a situation where it is not known in advance when the sleeping experts should be used. Therefore, we propose an online and adaptive method to wake up biased experts with Gradient Boosted Regression Trees (GBRT) \citep{chen_xgboost_2016} depending on a threshold. This enabled us to outperform the temperature predictions from Part I by making the EA more reactive, which was our main objective.

In this paper, for the sake of simplicity, we focus on the EA that achieved the best scores in Part I, namely Bernstein Online Aggregation (BOA) \citep{wintenberger_stochastic_2024} which is our benchmark. Furthermore, the best scores on average in Part I are without a sliding window, that is why we did not use sliding windows in this study, despite the fact that they can make EAs more reactive.

In Section \ref{sleeping_expert_framework}, we recall the SEF and explain how to adapt it to BOA. In Section \ref{how_we_wake_up_the_sleeping_experts}, we present the data and the specialized experts that we combined (the same experts as in Part I). Then we explain how we wake up the sleeping experts with GBRT and we provide a regret bound that depends on the performance of the GBRT. In Section \ref{results} we present our results and compare BOA with and without the SEF. We start with the tuning of the GBRT hyperparameters along with the tuning of the threshold used to wake up the specialized experts. Then we compare the reactivity and the weight behavior of BOA with and without the SEF. This is followed by the presentation of the main scores that we obtained and some differences between the stations and the lead times. Finally, we look at a score focusing on critical winter situations and an event where several large consecutive errors occurred without the SEF.

\section{Sleeping Expert Framework}
\label{sleeping_expert_framework}
\subsection{Changing environments}
\label{changing_environments}
In Part I, we used EA to outperform constant experts. \citet{devaine_forecasting_2013} tried to address a more challenging problem which is the framework of changing environments. In this framework, EA tries to take into account the fact that the best expert may change over time. The goal is no longer to beat the experts over the full sequence of observations but to beat them only on the iterations where they are good.

The key idea is that an expert can be specialized in certain events. So an expert could prove to be good at certain times for certain situations and bad the rest of the time.

The Specialized Framework tries to take advantage of this knowledge by using - by "waking up" - the specialized experts when they should be good. Otherwise, the EA puts the experts to sleep and does not use them.

\subsection{Prediction}
\label{Prediction}
Now, we will recall the SEF \citep{blum_empirical_1997,freund_using_1997,devaine_forecasting_2013} and some notations. To adapt any EA strategy to the SEF, one can use the "abstention trick" \citep{adamskiy_closer_2012,mourtada_efficient_2017} which we will describe in this section.

The goal of EA with the SEF\footnote{For "standard" EA without the SEF (all experts are always awake, $E_t=\{1,\ldots,N\}$ for all $t\geqslant 1$), see Section 2 of Part I.} is to predict with the help of subsets $E_t\subset\{1,\ldots,N\}$ of $N$ "awake" experts, $t=1,\ldots,T$, a sequence of observations $y_1,\ldots,y_T\in \mathcal{Y}$, where $\mathcal{Y}\subset{\mathbb{R}}$ is the outcome set. We will explain in Section \ref{how_we_wake_up_the_sleeping_experts}.\ref{how_we_used_gbrt}. how we choose the set of awake experts $E_t$.

At iteration $t$, expert $i$ makes the prediction $x_{i,t}^s$, $i=1,\ldots N$ and the aggregation predicts $\widehat{y}_t^s \in \widehat{\mathcal{Y}}$, with $\widehat{\mathcal{Y}} \subset{\mathbb{R}}$ the decision set. Here and in the following, exponent $s$ refers to the SEF in order to avoid confusions with $x_{i,t}$, $i=1,\ldots N$ and $\widehat{y}_t$ the expert's and aggregation's prediction without the SEF. We will for example denote BOA in the SEF by $\text{BOA}^s$.

In our study, we only use convex aggregation strategies. Hence, the aggregation's prediction is

\begin{equation}
    \widehat{y}_t^s=\sum_{i=1}^N w_{i,t}^s x_{i,t}^s=\frac{\sum_{i=1}^N v_{i,t}^s x_{i,t}^s}{\sum_{i=1}^N v_{i,t}^s},
\end{equation}
with $w_{i,t}^s=v_{i,t}^s/\sum_{k=1}^N v_{k,t}^s$, $v_{i,t}^s\in \mathbb{R}_+$ such that $\sum_{i=k}^N v_{k,t}^s>0$, and $(w_1^s,\ldots,w_N^s)\in \mathcal{X}$ where $\mathcal{X}=\{\mathbf{w}^s=(w_1^s,\ldots,w_N^s) | \sum_{i=1}^N w_{i}^s = 1, w_{i}^s \geqslant 0, i=1,\ldots,N\}$ is the set of convex weight vectors. $v_{i,t}$, $i=1,\ldots,N$, $t \geqslant 1$, are the experts' weights computed by the EA before normalization.

In the SEF, only a (non empty) set $E_t\subset \{1,\ldots,N\}$ of awake experts is used by the aggregation to make the prediction $\widehat{y}_{t}^s$. To do so, if expert $i\in E_t$ is awake we set
\begin{equation}
    x_{i,t}^s=x_{i,t}\label{eq 2}
\end{equation}

and if expert $j\notin E_t$ is sleeping (we tacitly assume that $\sum_{i\in E_t} v_{i,t}^s >0$ and $\sum_{i\in E_t} w_{i,t}^s >0$) we set

\begin{equation}
    x_{j,t}^s=\frac{\sum_{i\in E_t} v_{i,t}^s x_{i,t}}{\sum_{i\in E_t} v_{i,t}^s}.\label{eq 3}
\end{equation}

This is called the abstention trick \citep{mourtada_efficient_2017}. It ensures that the EA does not use at all\footnote{Note that the SEF is a special case of aggregation with experts that report their confidence which was introduced by \citet{blum_external_2007}. It has been adapted to aggregations with second-order regret bounds by \citet{gaillard_second-order_2014}. In \citet{gaillard_second-order_2014} an expert can be awake and completely used (confidence one), asleep and not used at all (confidence zero), or partially awake and partially used (confidence between zero and one). In our study, we do not use partially awake experts, only fully awake or asleep experts.} the sleeping experts for the prediction at iteration $t$ because in this case we have for $t \geqslant 1$

\begin{equation}
\label{eq_abstention_trick_1}
    \widehat{y}_{t}^s=\frac{\sum_{i\in E_t} v_{i,t}^s x_{i,t}}{\sum_{i\in E_t} v_{i,t}^s}
    =\frac{\sum_{i\in E_t} w_{i,t}^s x_{i,t}}{\sum_{i\in E_t} w_{i,t}^s},
\end{equation}

and

\begin{equation}
\label{eq_abstention_trick_2}
    x_{j,t}^s=\widehat{y}_{t}^s \text{ if } j\notin E_t,
\end{equation}

with $\widehat{y}_{t}^s$ the prediction of the EA, $E_t$ the set of awake experts, $w_{i,t}^s$ the weight of expert $i$ and $x_{i,t}^s$ its prediction at iteration $t$. That is to say that, by \eqref{eq_abstention_trick_2}, the sleeping experts are to yield the same prediction as the EA and that, by \eqref{eq_abstention_trick_1}, they have no influence on the prediction of the aggregation.

\subsection{Losses}
After the aggregation's prediction $\widehat{y}_t^s$, the observation $y_t$ is revealed and the aggregation incurs the loss $\ell(\widehat{y}_t^s,y_t)$, expert $i \in E_t$ incurs the loss $\ell(x_{i,t}^s,y_t)$, and expert $j\notin E_t$ incurs the loss $\ell(x_{j,t}^s,y_t)=\ell(\widehat{y}_t^s,y_t)$ where $\ell:\widehat{\mathcal{Y}} \times\mathcal{Y} \rightarrow \mathbb{R}$ is the loss function, a function which is convex and (sub)differentiable in its first argument. Thus, if expert $j\notin E_t$ is sleeping at iteration $t$, then, by equation \eqref{eq_abstention_trick_2}, it has the same loss as the aggregation at $t$. As in Part I, we will use the convex loss function $\ell(x,y)=(x-y)^2$ for each prediction $x \in \widehat{\mathcal{Y}}$ and observation $y \in \mathcal{Y}$.

Similarly to Part I, in the SEF, we denote $\ell_t:\mathcal{X} \rightarrow \mathbb{R}$  the function defined by $\ell_t(\mathbf{q})=\ell(\sum_{i=1}^N q_{i} x_{i,t}^s,y_t)$ where the $x_{i,t}^s$ are as in \eqref{eq 2} and \eqref{eq 3}. Hence, the loss of the EA is $\ell_t(\mathbf{w}_t^s)=\ell(\widehat{y}_t^s,y_t)$ and the loss of expert $i$ is $\ell_{t}(\mathbf{\delta}_i)=\ell(x_{i,t}^s,y_t)$ where $\mathbf{\delta}_i \in \mathbb{R}^N$ is the weight vector with 1 on coordinate $i$ and zero elsewhere.

Once the losses in the SEF are known, the aggregation can update $\mathbf{w}_t^s$ to $\mathbf{w}_{t+1}^s$ according to the EA strategy as explained in Part I. For example, to adapt BOA to the SEF, one only has to replace $x_{i,t}$ by $x_{i,t}^s$ for each iteration $t$ and expert $i=1,\ldots,N$ in Algorithm 3 of Part I. Algorithm \ref{specialized_expert_aggregation} summarizes the SEF.

Now, we define the cumulative loss of the aggregation in the SEF at iteration $t \geqslant 1$ as $\mathcal{L}_{t}^s=\sum_{k=1}^t \ell_{k}(\mathbf{w}_k^s)=\sum_{k=1}^t \ell(\sum_{i=1}^N w_{i,k}^s x_{i,k}^s,y_k)=\sum_{k=1}^t \ell(\sum_{i\in E_k} \frac{w_{i,k}^s}{\sum_{j \in E_k}w_{j,k}^s} x_{i,k},y_k)$ and the one of expert $i$ as $\mathcal{L}_{t}^s(\mathbf{\delta}_i)=\sum_{k=1}^t \ell_{k}(\mathbf{\delta}_i) \textbf{1}_{i\in E_k} + \ell_{k}(\mathbf{w}_k^s) \textbf{1}_{i\notin E_k}$ where $\textbf{1}$ is the indicator function (i.e. $\textbf{1}_{i\in E_k}$ equals $1$ if $i\in E_k$ and $0$ otherwise; similarly for $\textbf{1}_{i\notin E_k}$).

\subsection{Regret}
The regret of the aggregation in the SEF is denoted by $\mathcal{R}_{T}^s$, and as in Part I, the regret is the difference between the cumulative loss of the aggregation and the cumulative loss of a reference forecast (see Part I Sections 2.a. and 3.b. for more details).

\begin{algorithm}
\caption{Sleeping expert aggregation \citep{freund_using_1997}.}
\label{specialized_expert_aggregation}
\begin{algorithmic}
        \State Initialization: $v_{1,1},\ldots,v_{N,1}\in \mathbb{R}$
        \For{$t=1,\ldots,T$}
            \begin{enumerate}
                \item Choice of $E_t\subset\{1,\ldots,N\}$, the set of awake experts.
                \item Awake expert $i\in E_t$ makes the prediction $x_{i,t}^s=x_{i,t}$. Sleeping expert $j \notin E_t$ makes the prediction $x_{j,t}^s=\frac{\sum_{i\in E_t} v_{i,t}^s x_{i,t}^s}{\sum_{k\in E_t} v_{k,t}^s}$.
                \item The aggregation makes the prediction $\widehat{y}_{t}=\frac{\sum_{i\in E_t} v_{i,t}^s x_{i,t}^s}{\sum_{k\in E_t} v_{k,t}^s}$.
                \item The environment reveals the true outcome $y_{t}$.
                \item The aggregation incurs the loss $\ell_{t}(\mathbf{w}_t^s)$ and expert $i$ the loss $\ell_{t}(\mathbf{\delta}_i)$ if $i \in E_t$ and $\ell_{t}(\mathbf{w}_t^s)$ if $i \notin E_t$ with $\mathbf{w}_t^s=\left(\frac{v_{1,t}^s}{\sum_{i=1}^N v_{i,t}^s},\ldots,\frac{v_{N,t}^s}{\sum_{i=1}^N v_{i,t}^s}\right)$.
                \item The EA updates $v_{i,t} \longrightarrow v_{i,t+1}$, $i=1,\ldots,N$.
            \end{enumerate}

        \EndFor

\end{algorithmic}
\end{algorithm}

So after $T \ge 1$ iterations, the regret against expert $i=1,\ldots,N$ is \citep{devaine_forecasting_2013}

\begin{align}
    \mathcal{R}^s_{T}(\mathbf{\delta}_i)=\mathcal{L}_{T}^s - \mathcal{L}_{T}^s(\mathbf{\delta}_i)=\sum_{t=1}^T (\ell_{t}(\mathbf{w}_t^s) - \ell_{t}(\mathbf{\delta}_i)){\textbf{1}}_{i \in E_t}
\end{align}

where $\textbf{1}$ is the indicator function, $\ell_t(\mathbf{w}_t^s)$ the EA loss, and $\ell_t(\mathbf{\delta}_i)$ the loss of expert $i$. And the regret against a fixed convex combination $\mathbf{q}=(q_1,\ldots,q_N) \in \mathcal{X}$ of experts is \citep{devaine_forecasting_2013}

 \begin{align}
     \mathcal{R}^s_{T}(\mathbf{q})=\sum_{t=1}^T (\ell_t(\mathbf{w}_t^s) - \ell_{t}(\mathbf{q}^{E_t}))q(E_t)
 \end{align}
with $\ell_t(\mathbf{w}_t^s)$ the loss of the EA and $\ell_{t}(\mathbf{q}^{E_t})$ the loss of the convex combination $\mathbf{q}\in \mathcal{X}$ of experts, and $q(E_t)=\sum_{i \in E_t} q_i$ and $\mathbf{q}^{E_t}=(0,\ldots,0)$ if $q(E_t)=0$, and $\mathbf{q}^{E_t}=(\frac{q_1 \textbf{1}_{1\in E_t}}{q(E_t)},\ldots,\frac{q_N \textbf{1}_{N\in E_t}}{q(E_t)})$ if $q(E_t) > 0$. If all experts are always awake then $q(E_t)=1$ for all $t \geqslant 1$, and we retrieve the definition of the regret presented in Part I.

Thanks to a generic reduction introduced by \citet{adamskiy_closer_2012} and generalized by \citet{gaillard_second-order_2014} (Section 4 and Remark 6), one can deduce regret bounds in the SEF from the standard setting without the SEF.  So, for the adaptive version of BOA presented in Part I Section 2.c., in the SEF, with $\ell$ a differentiable loss function in its first argument and the gradient trick\footnote{ We recall that if the loss function is convex and subdifferentiable in its first argument, then one can use the gradient trick by replacing the loss $\ell_t(\mathbf{q})$ by the linear loss $\nabla \ell_t(\mathbf{w}_t) \cdot \mathbf{q}$, with $\mathbf{q}=(q_1,\ldots,q_N) \in \mathcal{X}$ a convex weight vector, $\nabla \ell_t(\mathbf{w}_t)$ a subgradient of $\ell_t$ at $\mathbf{w}_t$, $y_t$ the observation, $\widehat{y}_t$ the prediction of the EA at $t$, and $\cdot$ the inner product in $\mathbb{R}^N$.} (\citet{cesa-bianchi_prediction_2006} Section 2.5.), the regret bound against the fixed convex combination of experts $\mathbf{q} \in \mathcal{X}$ becomes

\begin{align}
    \mathcal{R}^s_{T}(\mathbf{q}) \leqslant O \left( \sum_{i=1}^N q_i \sqrt{ \sum_{t=1}^T \textbf{1}_{i\in E_t} \left(  \frac{ \partial \ell(\widehat{y}_t^s,y_t)}{\partial \widehat{y}_t^s}(x_{i,t}^s-\widehat{y}_t^s)\right)^2 } \right)
\end{align}
where  $\widehat{y}_t^s$ and $x_{i,t}^s$ are respectively the predictions of the aggregation and expert $i$, $y_t$ is the observation and $E_t$ is the set of awake experts at iteration $t$. This is a simplified version of the bound, for more details cf. Theorem 4.1 of \citet{wintenberger_stochastic_2024}.

By defining a compound expert $i^T$ as in \citet{mourtada_efficient_2017} as a sequence of experts\footnote{The compound expert $i^T=(i_1,\ldots,i_T)$ predicts at iteration $t$ like expert $i_t\in \{1,\ldots,N\}$. A compound expert could be $i^5=(3,2,2,2,9)$, which predicts like expert 3 at $t=1$ then like expert 2 at iterations 2, 3 and 4 and like expert 9 at iteration 5. And $i^T=(i,\ldots,i)$ is equivalent to the constant expert $i\in \{1,\ldots,N\}$.} $i^T=(i_1,\ldots,i_T)\in \{1,\ldots,N\}^T$, one can link the regret for the loss $\ell$ in the SEF with the compound expert $i^T$
 \begin{align}\label{eq:regretsef}
    \mathcal{R}^{s}_{T}(i^T) & =\sum_{t=1}^T ( \ell_{t}(\mathbf{w}_{t}^s) - \ell_{t}(i^T) ) {\textbf{1}}_{i_t \in E_t},
\end{align}

and if the compound expert is always awake, then for all $t\geqslant 1$, $i_t \in E_t$, and
 \begin{align}\label{eq:regretchang}
    \mathcal{R}^{s}_{T}(i^T) & =\sum_{t=1}^T \ell_{t}(\mathbf{w}_{t}^s) - \ell_{t}(i^T),
\end{align}

where $\textbf{1}$ is the indicator function, $\ell_{t}(\mathbf{w}_{t}^s)$ the loss of the EA, and $\ell_{t}(i^T)=\ell_t(\delta_{i_t})$ the loss of $i^T$ at iteration $t$.

Notice that if for all $t\geqslant 1$, $E_t=\{i_t\}$, then the EA predicts like $i^T$. We therefore reformulate the regret in a changing environment \eqref{eq:regretchang} as the SEF regret \eqref{eq:regretsef} given that the best expert is always active: $i_t \in E_t$, $t\ge 1$. Hence, one wants to reduce $E_t\subset \{1,\ldots,N\}$ aiming at the best strategy $E_t=\{i_t\}$. However, this is almost impossible to achieve, as this amounts to knowing in advance which is the best expert for each iteration. This is why we simplified this framework as explained in Section \ref{how_we_wake_up_the_sleeping_experts}., where we also provide a regret bound depending on how we wake up the experts.

\section{How we wake up the sleeping experts}
\label{how_we_wake_up_the_sleeping_experts}

\subsection{Why be more reactive}

We have seen in Part I that the various EAs perform well on average, consistent with the regret bounds of Section 2 in Part I. But at some point, this strength also becomes a disadvantage.

Indeed, these EAs are built to be good on average by minimizing the cumulative loss and by trying to converge to an oracle, either the best expert or the best fixed convex combination of experts. Due to this fact and despite being online and adaptive, these aggregations have trouble adapting to rapid changes in the set $E_t$, the set of experts that should be good at iteration $t$.

We can see in Figure 6 of Part I that for all EAs, the biased experts quickly have weights close to zero and become negligible. Furthermore, on this figure, all the weights become more and more stable and the aggregations become less reactive.  By reactive, we mean the ability to have fast and large weight variations.

The low reactivity of the EAs with second-order regret bounds can be partly explained by the excess loss\footnote{We recall the definition from Part I: for $t \geqslant 1$, expert $i=1,\ldots,N$, the excess loss is $\ell_{i,t}^{exc}= \ell_{t}(\mathbf{\delta}_i) - \sum_{j=1}^N w_{j,t} \ell_{t}(\mathbf{\delta}_j)$} $\ell_t^{exc}$. Indeed, the excess loss tends to stabilize the weights of the aggregation by favoring the experts that are close to the aggregation's prediction. Furthermore, all the learning rates used in Part I favor the experts with low cumulative losses, preventing previously bad experts from quickly regaining an important weight.

Hence, all these aggregation strategies rapidly lose the possibility of adding a lot of weight to an expert which was bad in the past but which becomes good for a short period of time, such as the $10\%$ quantile of the PEARP in mid-December in Chamonix, as can be seen in Figure 3 of Part I.

A more ambitious goal in a changing environment such as in our case is to try to compete with the experts only when they are good, thanks to the SEF. As mentioned before, the idea of the SEF is that each expert could have a specialty, which means that they should be better for some kind of events, and that they should only be used for those events.

\subsection{Our set of specialized experts}
\subsubsection{Data}
Our experts are the same as in Part I: the raw and post-processed versions of Integrated Forecasting System (IFS), forecasts of Application of Research to Operations at Mesoscale (AROME), Action de Recherche Petite Echelle Grande Echelle (ARPEGE) and the 10\%, 30\%, 50\%, 70\% and 90\% quantiles of the post-processed PEARP. They are respectively denoted raw.ifs, ppm.ifs, raw.aro, ppm.aro, raw.arp, ppm.arp and Q10, Q30, Q50, Q70 and Q90. We did not have access to other quantiles, which could have improved the EA. 

\begin{figure}[h]
\centerline{\includegraphics[width=19pc]{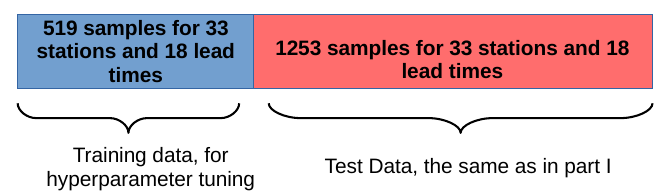}}
  \caption{Data used for this study.}\label{figure_training_data_papier}
\end{figure}

The data is split into a test set, containing the same data as in Part I, and a training set (see Figure \ref{figure_training_data_papier}). We recall that the test data (2020-03-30 to 2023-09-03) contains the predictions of the 00 TU run of all these experts and the corresponding temperature observations (obs) 2 meters above the ground for the 33 stations visible in Figure 1 of Part I. The predictions are for the lead times 6, 9, 12, 15, 18, 21, 24, 27, 30, 33, 36, 39, 42, 45, 48 for all experts plus the lead times 57, 72, 84 for the raw and post-processed ARPEGE, IFS and PEARP quantiles. The training data is the same, but from 2023-09-04 to 2025-02-03.

\subsubsection{Specialization}
In the SEF, we consider that our experts are specialized for different types of events. We consider the Q10 and Q30 experts to be specialized for events such as badly predicted cold spells. Or more generally, when the EA without the SEF makes too hot predictions and the aggregation's error is positive.

The Q70 and Q90 experts are considered to be specialized for events where the unbiased experts often make predictions which are too cold and the aggregation's error without the SEF is negative, such as in badly predicted heat waves.

The other experts (the Q50 expert and the post-processed and raw outputs of AROME ARPEGE and IFS) are considered to be specialized for daily forecasts and are always awake.

Hence, the goal of the biased experts is to increase or decrease the EA's prediction at the right time, when the unbiased experts are not enough for the EA to make an accurate prediction.

\subsection{How to be more reactive}
So, as in \citet{devaine_forecasting_2013} and \citet{vyugin_online_2019}, we introduce specialized/sleeping experts. There is a big difference however between these papers and our study. It is crucial to see that in \citet{devaine_forecasting_2013} and \citet{vyugin_online_2019}, they know in advance $E_{t+1}$, the set of experts which should be the best for the next iteration. Indeed, their experts are specialized for weekends, holidays, summer or winter for example. So it is easy to wake up or not wake up, as the case may be, the appropriate set of experts at the right time.

In our case, the difficulty stems from the fact that we do not know $E_{t+1}$ at iteration $t\geqslant1$. In other words, we do not know a priori when the aggregation will make predictions which are too hot or too cold. So, at iteration $t$, we will have to predict $\widehat{E}_{t+1}$, the set of experts that should be good for iteration $t+1$.

One solution would be to have a human forecaster who, thanks to his knowledge and understanding of the weather and the NWP models, would know when to wake up or put to sleep the different experts. In practice, however, this proves to be very difficult given the frequency of runs and the number of lead times and stations involved. For example, at Meteo France, 4 times a day, a specific aggregation is performed for each pair of more than 2500 stations and 32 lead times.

Another solution is to use machine learning techniques to estimate $E_t$ in order to know when to use the biased experts. It is not obvious which machine learning method to use; for example, one could use linear discriminant analysis \citep{zhao_linear_2024} or neural networks \citep{marzban_neural_2003,rasp_neural_2018}. 

We chose GBRT from the xgboost library of \citet{chen_xgboost_2016}. We chose this machine learning method in particular because it combines the advantages of tree-based and gradient boosted models. Unlike discriminant analysis, GBRT does not make assumptions on the distribution. Furthermore, its training is quite fast thanks to the built in parallelization. This is important since we want to train the GBRT before each prediction because we want an online algorithm. GBRT is also a well-known state-of-the-art machine learning algorithm that has produced good results in a wide variety of fields.

GBRT has also been widely used for weather prediction and analysis \citep{fei_factors_2020,potdar_toward_2021,mony_evaluating_2021,flora_using_2021,silva_using_2022} and for soil temperature prediction \citep{liu_correction_2022,nanda_soil_2020} outperforming various machine learning methods such as random forests among others. Finally, it has also been used for air temperature forecasts by \citet{ma_prediction_2020} but only for very short forecast horizons of less than three hours.

Several ways exist to deal with the sleeping experts, and there is no obvious choice. In our study, we chose to combine EA with GBRT in order to wake up the specialized experts and to make the EA more reactive.

\subsection{How we used GBRT}
\label{how_we_used_gbrt}

At each iteration $t \geqslant 100$, for each station-lead time couple, we trained a GBRT in order to predict $\widehat{e}_t$ which tries to estimate $e_t=\widehat{y}_t - y_t$ the aggregation's error of BOA without the sleeping experts Q10, Q30, Q70 and Q90. For each station–lead time pair, the GBRT was trained using all historical data available up to time $t$, excluding data from other lead times and from different stations.

We did not use the GBRT and put the specialized experts to sleep for $t < 100$ (arbitrary choice), because we considered that fewer than 100 training samples are not enough to train a GBRT.

We also oversampled 2 times (arbitrary choice) all the observations with an absolute error larger than a threshold $e^\text{thr}$ (we explain how we chose $e^\text{thr}$ in section \ref{results}.\ref{hyperparameter}.). This distorts the training distribution, so that the GBRT model learns better the samples with large errors, which are difficult to predict.

In Appendix A, there is a complete list of the variables used to train the GBRT. To sum up, we tried to choose variables that correspond to what a human forecaster would do at the beginning of the working day: mainly look at the correlation between weather forecasting models and observations, and see if the models agree for a scenario.

There are therefore mainly different spreads and variances between the experts and also the difference between the experts and the observed temperature at the first lead time. This is possible and makes sense since in practice, the predictions of the raw and post-processed NWP models are available after the first lead time has been observed. Furthermore, in this study we only assess the EA after the second lead time which is 6 hours.

\begin{figure*}[h]
 \centerline{\includegraphics[width=39pc]{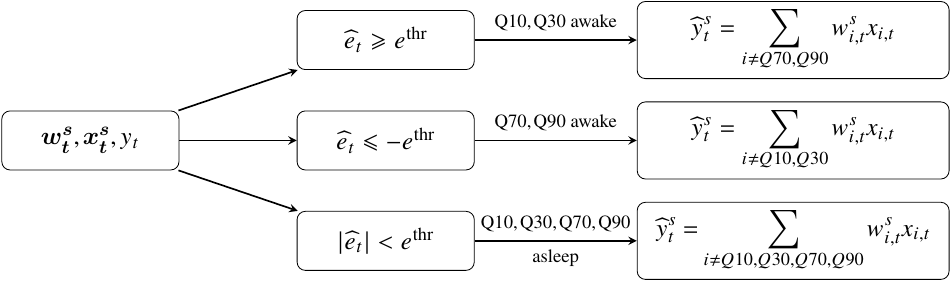}}
  \caption{Scheme illustrating the activation rule for specialized experts based on the gradient boosted regression tree's estimate $\widehat{e}_t$ of the aggregation error and the threshold $e^\text{thr}$: if $\widehat{e}_t \geqslant e^\text{thr} $, experts Q10 and Q30 are woken up; if $\widehat{e}_t \leqslant -e^\text{thr} $, experts Q70 and Q90 are woken up; otherwise, all biased experts remain asleep.}\label{shema_biased}
\end{figure*}

Finally, to choose $E_t$ and wake up the specialized experts or not, we used a simple rule depending on the threshold $e^\text{thr}$, illustrated in Figure \ref{shema_biased}: if $\widehat{e}_t \leqslant -e^\text{thr}$ we wake up Q70 and Q90, if $\widehat{e}_t \geqslant e^\text{thr}$ we wake up Q30 and Q10 and if $|\widehat{e}_t| < e^\text{thr}$ all the biased experts are asleep. The Q50, and the post-processed and raw outputs of AROME, ARPEGE and IFS are always awake. Thanks to this rule, we managed to run $\text{BOA}^{s}$ in the SEF in a fully online adaptive way.

In our SEF, the specialized experts combined with the GBRT try to predict the extreme events, whilst the other experts try to predict the average behavior of the temperature. The success of the SEF is strongly related to whether or not the specialized experts are indeed specialized for the events they are awakened for.

So, if the GBRT works well, at $t\geqslant 100$, when the biased experts are not woken up, the loss of the EA is bounded by $\ell_t \leqslant ({e^{\text{thr}}})^2$. And when the specialized experts are woken up, if they are "specialized enough " for iteration $t$, they will improve the forecast. 

We can see however that there is positive feedback between the GBRT and the aggregation. If the GBRT often wakes up biased experts at iterations where they increase the error of the aggregation, the aggregation would decrease the weight of these biased experts, and vice versa. So, as with all experts, the weights of the specialized experts will tend towards an optimum adapted to the quality of the GBRT, regardless of the predictions and observations. 

\subsection{Regret bound depending on the GBRT classification}
\label{regret_bound_classification}

As explained previously, to simplify the problem of predicting $E_t$, we do not try to predict the error $e_t=\widehat{y}_t - y_t \in \mathbb{R}$ where $\widehat{y}_t$ is the prediction of BOA without the biased experts at iteration $t \geqslant 1$. Instead, we only try to predict $k=3$ cases with GBRT (it could be any other algorithm): GBRT tries to predict $e_t' \in \{1,2,3\}$, where $e_t'=1$ if $e_t \leqslant -e^\text{thr}$, $e_t'=2$ if $|e_t| < e^\text{thr}$ and $e_t'=3$ if $e_t \geqslant e^\text{thr}$.

Let GBRT predict $\widehat{e_t}' \in \{1,2,3\}$, and  $\widehat{E}_t \subset \{ 1,\ldots,N \}$ be the set of experts that the GBRT wakes up at $t$, with $\widehat{E}_t=\widehat{E}^k$, where $\widehat{E}^k$ is the set of experts that is woken up when $\widehat{e_t}'=k$, $k=1,2,3$. For example, in our case, $\widehat{E}^1=\{\text{Q70,Q90, raw.aro, raw.arp, raw.ifs, ppm.aro, ppm.arp, ppm.ifs, Q50}\}$.

Since we always want to make a temperature prediction, it makes sense to suppose that there is always an awake expert (for all $t\geqslant 1$, $E_t$ and $\widehat{E}_t$ are non empty) and that $i^T=(i_1,\ldots,i_T)$ is a compound expert such that for all $t\geqslant 1, i_t \in E_t$. In this case

\begin{align*}
    \mathcal{R}^s_{T}(i^T) =& \sum_{t=1}^T ( \ell_{t}(\mathbf{w}^s_t) - \ell_{t}(\mathbf{\delta}_{i_t}) ) \textbf{1}_{i_t\in E_t}\\
     =& \sum_{t=1}^T \ell_{t}(\mathbf{w}^s_t) - \ell_{t}(\mathbf{\delta}_{i_t}) \\
      =& \sum_{t=1}^T \sum_{\widehat{k}=1}^3 \left[ \ell \left( \sum_{i\in\widehat{E}^{\widehat{k}}} w_{i,t}^s x_{i,t}, y_t \right) - \ell_{t}(\mathbf{\delta}_{i_t}) \right] \textbf{1}_{\widehat{e_t}' = \widehat{k}}\\
      \leqslant& \sum_{t=1}^T \sum_{{\widehat{k}}=1}^3 \left[ \sum_{k=1}^3 \left( \max_{\widehat{e_t}'=\widehat{k},e_t'=k} \ell_{t}(\mathbf{w}^s_t)  - \ell_{t}(\mathbf{\delta}_{i_t})  \right) \textbf{1}_{e_t'=k} \right] \textbf{1}_{\widehat{e_t}' = \widehat{k}}\\
      \leqslant& \sum_{t=1}^T \sum_{{\widehat{k}}=1}^3 \left[ \sum_{k=1}^3 \left( \max_{\widehat{e_t}'=\widehat{k},e_t'=k} r_{i_t,t}(\mathbf{w}^s_t) \textbf{1}_{e_t'= k} \right) \right] \textbf{1}_{\widehat{e_t}' = \widehat{k}},
\end{align*}

hence
\begin{equation}
     \mathcal{R}^s_{T}(i^T) \leqslant \sum_{{\widehat{k}}=1}^3 \sum_{k=1}^3 n_{\widehat{k},k}^T \max_{\widehat{e_t}'=\widehat{k},e_t'=k} r_{i_t,t}(\mathbf{w}^s_t)
\end{equation}

where for $t \geqslant 1$, $r_{i_t,t}(\mathbf{w}^s_t)= \ell_{t}(\mathbf{w}^s_t) - \ell_t(\mathbf{\delta}_{i_t})$ the instantaneous regret, and $n_{\widehat{k},k}^t =\sum_{r=1}^t \textbf{1}_{\widehat{e}_r'=\widehat{k},e_r'=k}$, with $\widehat{k},k=1,2,3$. If the specialized experts are adapted to the situations for which they are woken up such that for $\widehat{k},k=1,2,3$, $\max_{\widehat{e}_t'=e_t'=k} r_{i_t,t}(\mathbf{w}^s_t) \leqslant \max_{\widehat{e}_t'=\widehat{k}, e_t'=k} r_{i_t,t}(\mathbf{w}^s_t) $, then one wants $n_{\widehat{k},k}^T=0$ for $\widehat{k} \ne k$ to minimize the regret bound. This indicates that an increasing number of false classifications ($n_{\widehat{k},k}^T$, $\widehat{k} \ne k$) leads to a higher regret bound.

If the predictions of the GBRT are perfect (for all $t \geqslant 1$, $\widehat{e}_t'=e_t'$ hence $n_{\widehat{k},k}^T=0$ for $\widehat{k} \ne k$), then
\begin{align}
    \mathcal{R}^s_{T}(i^T) \leqslant& \sum_{t=1}^T \sum_{\widehat{k}=1}^3 \left( \max_{\widehat{e}_t'=k} \ell_t(\mathbf{w}^s_t) - \ell_{t}(\mathbf{\delta}_{i_t}) \right) \textbf{1}_{\widehat{e}_t'=\widehat{k}}\\
    \leqslant& \sum_{{\widehat{k}}=1}^3 n_{\widehat{k}}^T \max_{\widehat{e}_t'=\widehat{k}} r_{i_t,t}(\mathbf{w}^s_t),
\end{align}

with $n_{\widehat{k}}^t=\sum_{r=1}^t \textbf{1}_{\widehat{e}_r'=\widehat{k}}$, for $t \geqslant 1$ and $\widehat{k}=1,2,3$.

One can easily generalize these results for $K\in \mathbb{N}$ categories and to any model predicting $\widehat{e}_t' \in \{1,\ldots,K\}$ for $t=1,\ldots,T$, by replacing the three cases by $K$ cases.

\section{Results}
\label{results}
In this part, we will on the one hand look at the results from a global point of view (all the stations and lead times) and on the other hand look at special cases like the Chamonix station for a 48-hour lead time. For the following, as in Part I, we used R 3.6.1 and the Opera 1.1.1  package \citep{gaillard_opera_2016} that we modified in order to deal with the sleeping  experts with the xgboost 1.7.7.1 package \citep{chen_xgboost_2016}.

\subsection{Threshold and hyperparameter tuning}
\label{hyperparameter}
A downside of making the aggregations more reactive could be to add some noise and degrade the predictions, especially by waking up the wrong sleeping expert at the wrong time. For example, waking up the Q90 (Q10) expert when the loss of the aggregation is positive (negative) is likely to make the prediction of the aggregation worse. That is why the choice of the threshold $e^\text{thr}$ is very important.

Hence, before using the GBRT, we had to choose a set of hyperparameters among the 15876 possible combinations of hyperparameters presented in Table \ref{table_hyperparameters} and a threshold $e^{\text{thr}}$ among $\{0, 0.5, 1, 1.5, 2, 2.5, 3, 3.5 \}$.

\begin{table*}[h]
\begin{center}
\begin{tabular}{cc}
\hline
Hyperparameters & Values \\
\hline
No. of rounds & 1, 2, 3, 4, 5, 6, \textbf{7}\\
Max depth & 4, 5, 6, 7, \textbf{8}, 9, 10\\
Learning rate & 0.2, 0.3, \textbf{0.4}, 0.5, 0.6, 0.7, 0.8, 0.9, 1\\
Min child weight & 5, 10, 20, 35, 60, \textbf{100}\\
\renewcommand{\arraystretch}{1}\begin{tabular}{@{}c@{}}Ratio of predictors randomly \\ selected per node\end{tabular} & 0.65, 0.7, 0.75, 0.8, \textbf{0.85}, 0.9\\
\hline
\end{tabular}
\caption{Hyperparameters used for the tuning of the gradient boosted regression trees. The selected hyperparameters for the threshold $e^\text{thr}=0.5$°C are in bold.}
\label{table_hyperparameters}
\end{center}
\end{table*}

Since we trained one aggregation for each station-lead time pair, at first glance the most natural would have been to tune one set of xgboost hyperparameters for each station-lead time couple. But we tuned only one set of hyperparameters for all the stations and lead times for two main reasons.

First of all, we wanted to compare our results to those of Part I. This implies that we only can tune the hyperparameters with a different data from Part I. For this reason we had to split our data as shown in Figure \ref{figure_training_data_papier}, leaving us with only 519 samples for each station-lead time pair to tune the xgboost hyperparameters. This turns out to be too few samples for tuning a separate set of hyperparameters for each station–lead time pair. 

Secondly, we wanted our work to be usable in an operational context, where there are many more stations (e.g. more than 2500 at Meteo France) and more than 18 lead times. This makes it really complicated to tune xgboost hyperparameters for every possible couple of stations and lead times, which could also lead to overfitting. Furthermore, it turned out that many sets of hyperparameters led to almost the same performance on the training data. Hence, xgboost does not appear to be highly sensitive to the choice of hyperparameters in our case.

That is why we chose for all the stations and lead times the same set of hyperparameters among those of Table \ref{table_hyperparameters} and the same threshold $e^{\text{thr}}$.

For each threshold $e^{\text{thr}} \in \{0, 0.5, 1, 1.5, 2, 2.5, 3, 3.5 \}$, we did a cross validation on the training data for every GBRT associated to the 15876 different sets of hyperparameters. We evaluated the GBRT with an Equitable\footnote{Here, equitable means that a model which always forecasts the same class has the same score as a model making random predictions.} Skill Score (ESS) \citep{gerrity_note_1992} which is recommended by the World Meteorological Organization \citep{wmo_manual_2010} for classification evaluation. For more details and to find out how to compute the ESS, see Appendix B. We denote $\text{GBRT}^{e^\text{thr}}$ the GBRT with the hyperparameters tuned for $e^{\text{thr}}$.

\begin{figure*} [h]
 \centerline{\includegraphics[width=19pc]{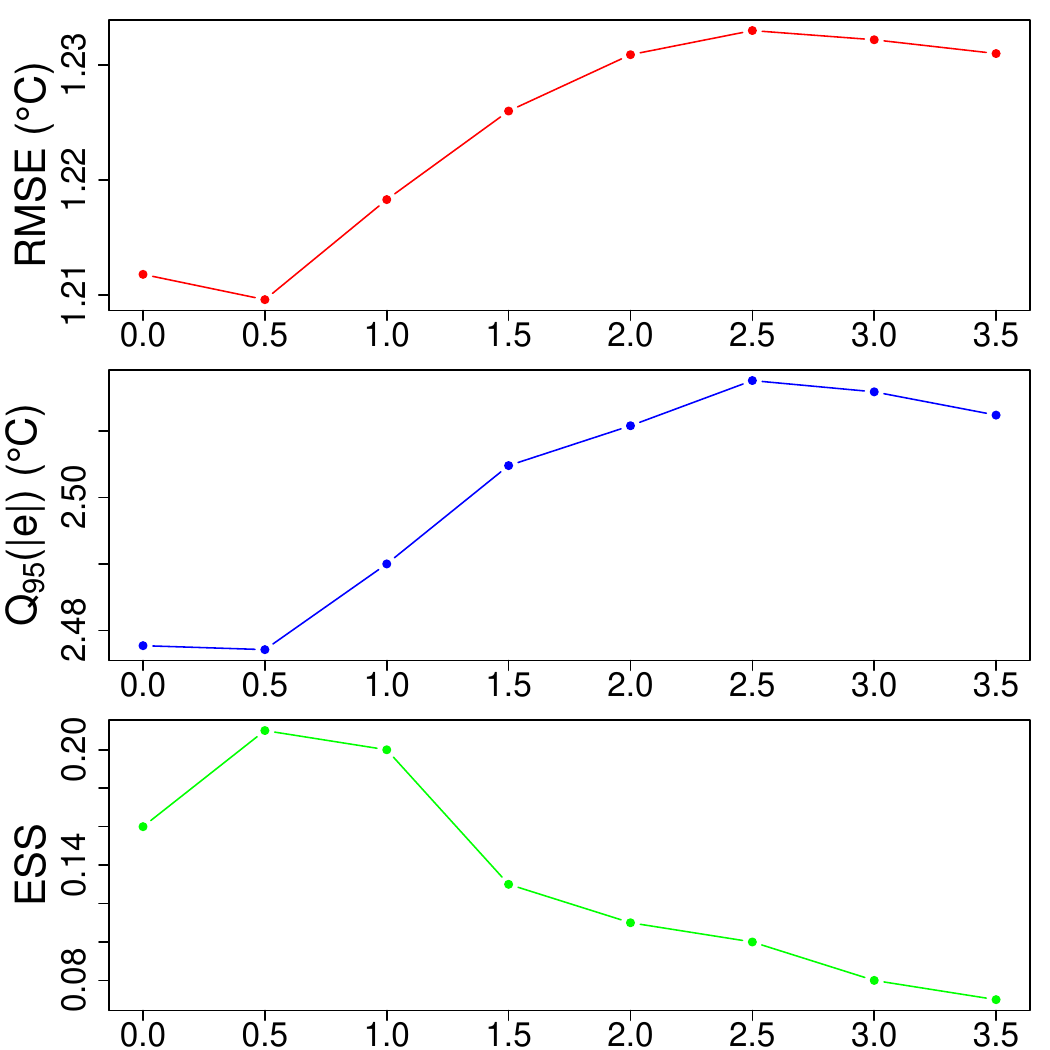}}
  \caption{Depending on the threshold $e^{thr} \in \{0, 0.5, 1, 1.5, 2, 2.5, 3, 3.5 \}$ : (bottom) Equitable Skill Score (ESS) of Gradient Boosted Regression Trees tuned for $e^{thr}$ ($\text{GBRT}^{e^{\text{thr}}}$), (middle) the 95\% quantile of the absolute error $Q95(|e|)$ of the EA with $\text{GBRT}^{e^{\text{thr}}}$ and $e^{thr}$ to wake up the sleeping experts, (top) the Root Mean Squared Error (RMSE) of the expert aggregation with $\text{GBRT}^{e^{\text{thr}}}$ and $e^{thr}$ to wake up the sleeping experts.}
  \label{compare_threshold_519}
\end{figure*}

Then, for each threshold $e^{\text{thr}} \in \{0, 0.5, 1, 1.5, 2, 2.5, 3, 3.5 \}$, for each station-lead time pair, we run in the SEF on the training data, an EA combined with $\text{GBRT}^{e^\text{thr}}$. The scores\footnote{For more details on the RMSE and the $Q95(|e|)$, see Part I, Section 3.c.} are presented in Figure \ref{compare_threshold_519}, where one can see the ESS of $\text{GBRT}^{e^\text{thr}}$, the Root Mean Squared Error (RMSE) and the 95\% quantile of the absolute error $Q95(|e|)$ of the EAs in the SEF combined with $\text{GBRT}^{e^\text{thr}}$ for $e^{\text{thr}} \in \{0, 0.5, 1, 1.5, 2, 2.5, 3, 3.5 \}$.

On the training data, we obtained the best results\footnote{We take this opportunity to mention that we tried to pool data together in order to improve the predictions. First we combined all data from all stations and lead times to train the GBRT, then we gathered all data of all the lead times of the considered station, and added the lead time and the station as predictors to train the GBRT. But it turned out that it added noise to the EA when used on the training data.} for all three scores ($\text{ESS}=0.21$, $\text{RMSE}=1.21$°C and $Q95(|e|)=2.48$°C) with the threshold $e^\text{thr}=0.5$°C and $\text{GBRT}^{0.5}$. So, in the following, we will test $\text{BOA}^s$ with $e^{\text{thr}}=0.5$°C and $\text{GBRT}^{0.5}$ (hereafter denoted as GBRT) on the test data.

\subsection{Weights behavior}
Without the SEF, we saw in Part I that the EAs were unable to quickly activate the biased experts that had too often performed poorly. Hence, the first goal of using the SEF is to make the weights of the biased experts (Q10, Q30, Q70 and Q90) more reactive than in Part I.

\begin{figure} [h]
\centering
    \includegraphics[width=33pc]{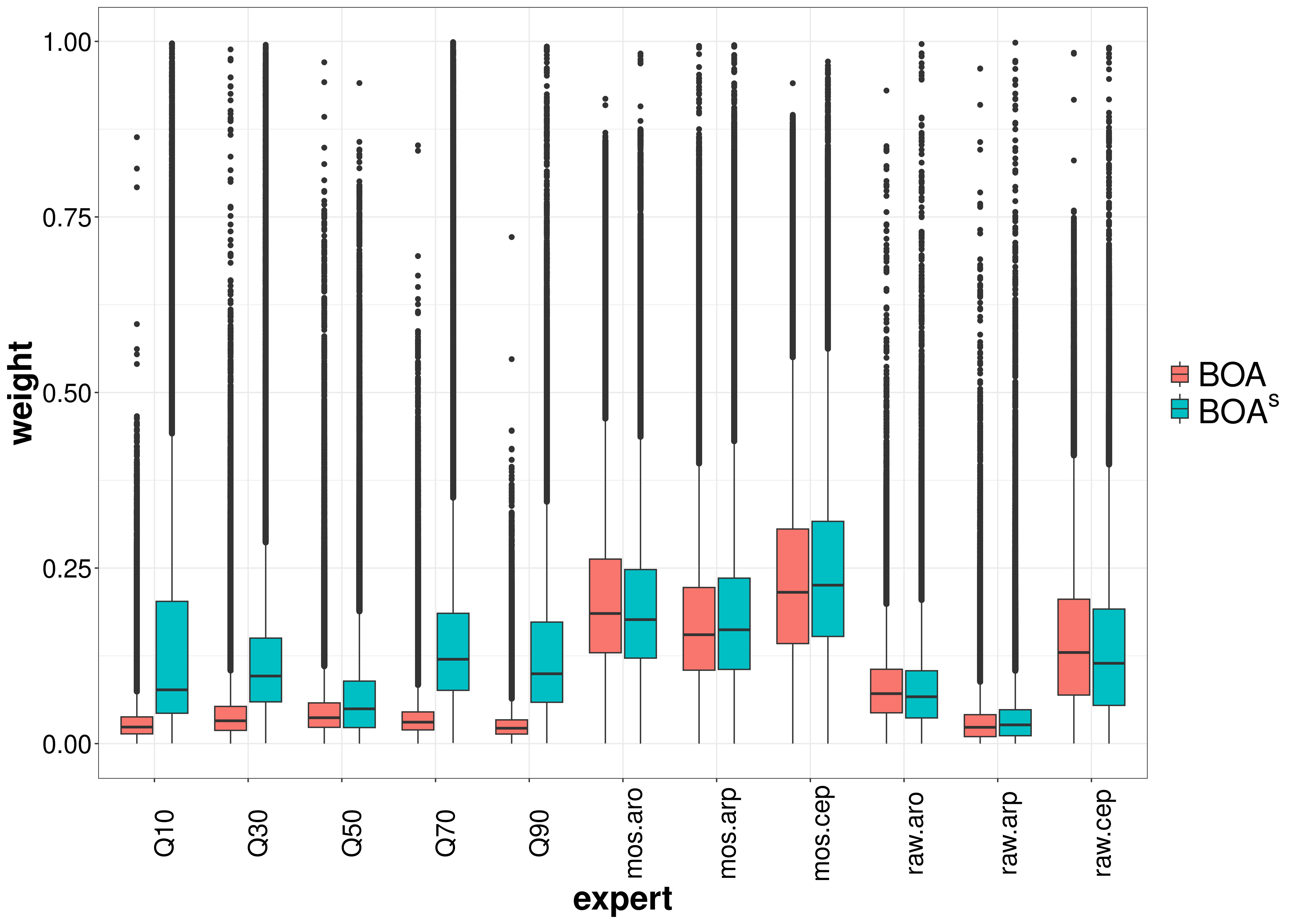}
\caption{Boxplot of the weight of the experts (when they are awake), for all the stations, lead times and iterations between 2020-03-30 and 2023-09-03, for BOA \citep{wintenberger_stochastic_2024} and BOA with the sleeping expert framework ($\text{BOA}^s$).}
\label{boxplot_weights_BOA_vs_sleep_BOA}
\end{figure}

In Figure \ref{boxplot_weights_BOA_vs_sleep_BOA}, one can see the boxplots of the weights of the experts for BOA and $\text{BOA}^s$ when they are awake, across all the stations, lead times and iterations. It is clear that overall, in the SEF, when the specialized experts are awake, they have much larger weights than without the SEF.

The weight distribution of the other experts does not change much between BOA and $\text{BOA}^s$. The small changes are probably mainly due to the fact that in the SEF, the weights of the sleeping experts are dispatched between the non sleeping experts.

\begin{figure*} [h]
 \centerline{\includegraphics[width=39pc]{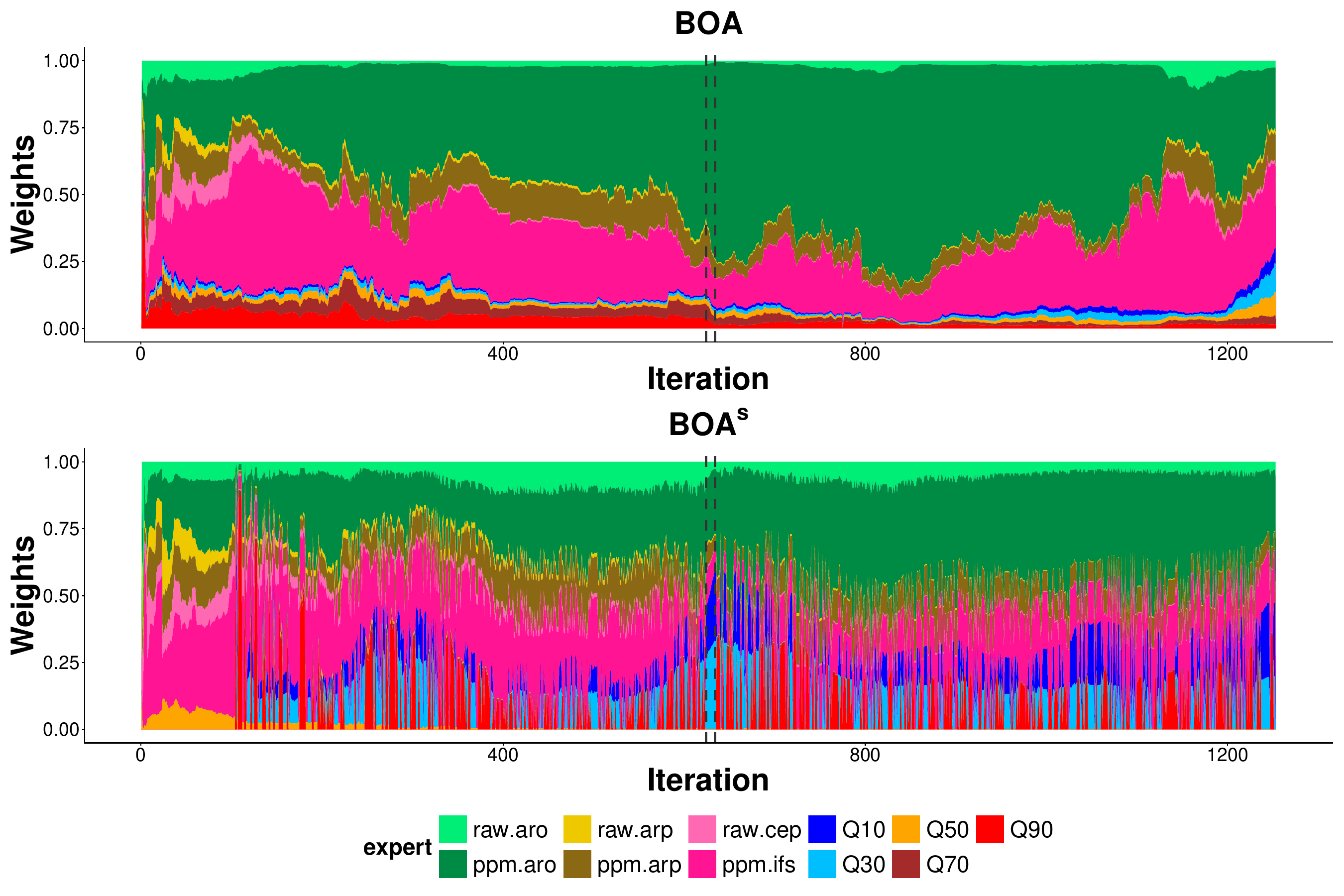}}
  \caption{Weights of the experts, for (top) BOA \citep{wintenberger_stochastic_2024}, and (bottom) BOA with the sleeping expert framework ($\text{BOA}^{s}$), in Chamonix, from 2020-03-30 to 2023-09-03, for a 48-hour lead time. The dashed grey vertical lines delimit the December period where all experts make too warm predictions except the Q10 and Q30 experts.}\label{BOA_sleep_BOA_window1253_ech48_74056001_2020-03-30_2023-09-03_weights}
\end{figure*}

The weights of the experts of BOA and $\text{BOA}^{s}$ can be seen in Figure \ref{BOA_sleep_BOA_window1253_ech48_74056001_2020-03-30_2023-09-03_weights} for a 48-hour lead time in Chamonix from March 2020 to September 2023 as in Figure 6 of Part I.

We can see in this example of Chamonix, that the SEF enables much faster and larger changes in the weights of the biased experts. Indeed, there are many abrupt changes of weights when the sleeping experts are woken up by the GBRT. This is visible in the bottom plot of Figure \ref{BOA_sleep_BOA_window1253_ech48_74056001_2020-03-30_2023-09-03_weights} where there are many blue and red vertical strokes, with even wider ones just around iteration 630 (between the dashed gray lines), corresponding to December 2021. We will take a closer look at this event in Section \ref{results}.\ref{special_case_study}.\ref{consecutive_errors}.).

The dark and bright blue strokes around iteration 630, show that $\text{BOA}^{s}$ puts large weights on Q10 and Q30 at these iterations where BOA makes several predictions which are too hot as we will see later.

So the SEF combined with GBRT, seems to allow $\text{BOA}^{s}$ to use the biased experts efficiently  even after a large number of iterations where the sleeping experts would have made bad predictions if they had been awake. 

\subsection{Scores}
With our set of hyperparameters and the threshold $e^{\text{thr}}=0.5$°C, the ESS across all the stations and lead times for the 1253 iterations between 2020-03-30 and 2023-09-03 is $ESS=0.23$.

We computed the RMSE and the 95\% quantile of the absolute error $Q95(|e|)$. The scores are presented in Table \ref{table_rmse_Q95}. We tested if the score differences are statistically significant with a Diebold–Mariano (DM) test \citep{diebold_comparing_1995} for the RMSE and a quantile test \citep{johnson_two-sample_1987} for the $Q95(|e|)$, both with a Benjamini-Hochberg (BH) procedure \citep{benjamini_controlling_1995}. The results of these tests can be seen in Tables \ref{table_diebold_mariano_rmse} and \ref{table_quantile_test_Q95}.

\begin{table*}[h]
\begin{center}
\begin{tabular}{cccc}
\hline
& RMSE (°C) & $Q95(|e|)$ \\
\hline
BOA & 1.24 & 2.53 \\
$\text{BOA}^s$ & 1.21 & 2.46 \\
\renewcommand{\arraystretch}{1}\begin{tabular}{@{}c@{}}Best convex\\Combination\end{tabular} & 1.23 & 2.50  \\
\hline
\end{tabular}
\caption{RMSE and $Q95(|e|)$ of BOA, $\text{BOA}^s$ and the oracle of the best fixed convex combination of experts.}
\label{table_rmse_Q95}
\end{center}
\end{table*}

\begin{table*}[h]
\begin{center}
\begin{tabular}{cccc}
\hline
& BOA & $\text{BOA}^s$ & \shortstack{Best convex\\Combination} \\
\hline
BOA &  & 0.0 & 0.0 \\
$\text{BOA}^s$ & 36.4 &  & 31.3 \\
\renewcommand{\arraystretch}{1}\begin{tabular}{@{}c@{}}Best convex\\Combination\end{tabular} & 78.1 & 23.1 &  \\
\hline
\end{tabular}
\caption{Percentage of station-lead time couples with a statistically significant difference between the RMSE between BOA, $\text{BOA}^s$ and the oracle of the best fixed convex combination of experts. On row $i$ and column $j$ is the percentage of lead time/station couples for which the one sided DM test rejects the null hypothesis in favor of the aggregation of row $i$ against the aggregation of column $j$ with a Benjamini–Hochberg (BH) procedure to take account of multiple testing with a 0.05 level. The tests were run with the dm.test function of the R package forecast 8.24.0 and the p.adjust function of the R package stats 3.6.1.}
\label{table_diebold_mariano_rmse}
\end{center}
\end{table*}

\begin{table*}[h]
\begin{center}
\begin{tabular}{cccc}
\hline
& BOA & $\text{BOA}^s$ & \shortstack{Best convex\\Combination} \\
\hline
BOA &  & 0 & 0 \\
$\text{BOA}^s$ & 5.9 &  & 5.1 \\
\renewcommand{\arraystretch}{1}\begin{tabular}{@{}c@{}}Best convex\\Combination\end{tabular} & 0 & 0 &  \\
\hline
\end{tabular}
\caption{Percentage of station-lead time couples with a statistically significant difference between the $Q_{95}(|e|)$ between BOA, $\text{BOA}^s$ and the oracle of the best fixed convex combination of experts. On row $i$ and column $j$ is the percentage of lead time/station couples for which the quantile test rejects the null hypothesis in favor of the aggregation of row $i$ against the aggregation of column $j$ with a Benjamini–Hochberg procedure to take account of multiple testing with a 0.05 level. The tests were run with the quantileTest function of the R package EnvStats 3.1.0 and the p.adjust function of the R package stats 3.6.1.}
\label{table_quantile_test_Q95}
\end{center}
\end{table*}

The RMSE of $\text{BOA}^{s}$ across all the stations and lead times is 1.21°C, whereas it is 1.24 for BOA without the SEF. So we could make the aggregation more reactive whilst improving the RMSE. Indeed, the improvement is statistically significant for 36.4\% of the station lead time couples.

Furthermore, $\text{BOA}^{s}$ outperformed the best fixed convex combination of experts oracle - that could not be outperformed in Part I - for 31.3\% of the station lead time couples (see Table \ref{table_diebold_mariano_rmse}).

Another important score is the 95\% quantile of the absolute error $Q_{95}(|e|)$, which is much more relevant with regards to extreme events and large errors. The $Q_{95}(|e|)$ across all the stations and lead times of BOA is 2.53°C whereas the one of $\text{BOA}^{s}$ is 2.46°C. The $Q_{95}(|e|)$ improvement thanks to the SEF is statistically significant for almost 6\% of the station lead time couples. Therefore, overall, using the SEF for our study helped to achieve fewer large errors.

The fact that the RMSE and the $Q_{95}(|e|)$ of $\text{BOA}^{s}$ across the 1253 observations 33 stations and 18 lead times is lower than it is for BOA, means that when the GBRT woke up some sleeping experts, it most often concerned the right experts at the right moment.

\begin{figure*} [h]
 \centerline{\includegraphics[width=39pc]{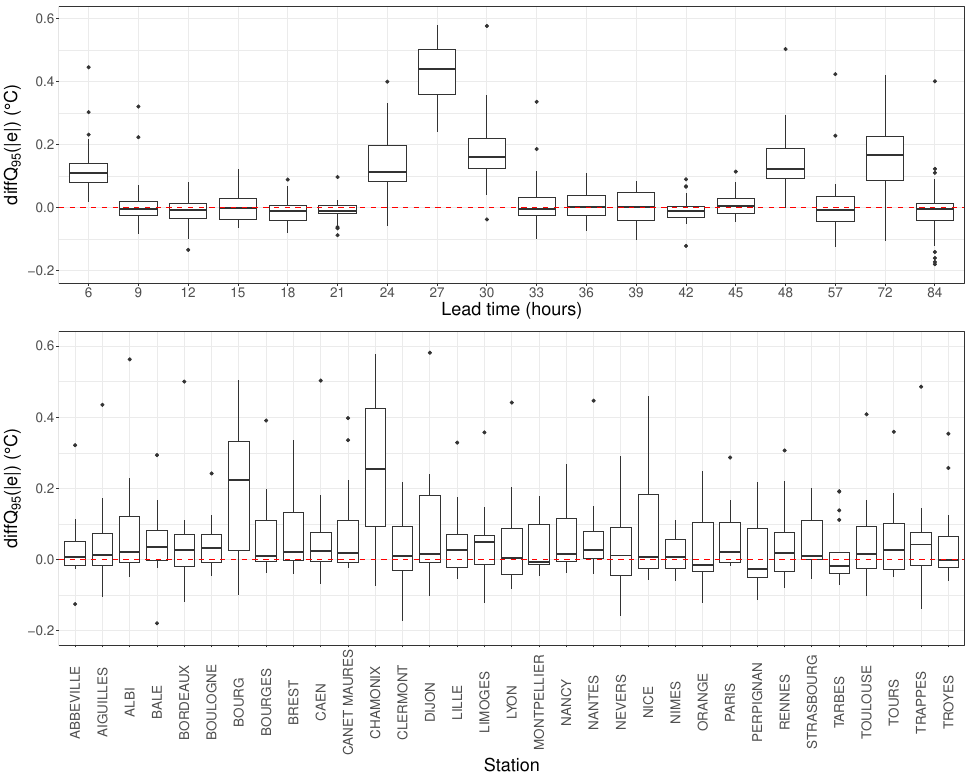}}
  \caption{Box plot of diff$Q_{95}(|e|)$ the difference between the 95\% quantile of the absolute error of BOA and BOA in the Sleeping expert framework ($\text{BOA}^{s}$). (Top) Depending on the lead time. (Bottom) Depending on the station. When the difference is positive (over the dashed red line), then $\text{BOA}^{s}$ is better. }\label{diffQ95_BOA_vs_BOAsleep}
\end{figure*}

In Figure \ref{diffQ95_BOA_vs_BOAsleep}, one can see the boxplots of the difference between the $Q_{95}(|e|)$ of BOA and $\text{BOA}^{s}$. So when the difference is positive, then $\text{BOA}^{s}$ has a better $Q_{95}(|e|)$ than BOA. In the upper plot of Figure \ref{diffQ95_BOA_vs_BOAsleep}, one can see that there is a periodicity in the contribution of the SEF. The improvement for night prediction is much higher than for day predictions. The 27 hour lead time is especially noteworthy, because it is by far the best lead time regarding the improvement due to the SEF. The 84 hour  lead time appears to be the worst.

The fact that the SEF with the GBRT does not always improve BOA relates to our remark that it is difficult to refine PPMs, particularly when they are aggregated, as their remaining errors are already small and nearly random.

The lower plot of Figure \ref{diffQ95_BOA_vs_BOAsleep} also shows that there is a variability of the SEF contribution depending on the stations.

For the large majority of the stations, the SEF improves the $Q_{95}(|e|)$ because the $Q_{95}(|e|)$ differences are most of the time positive as shown by the distributions. For Chamonix and Bourg (Saint Maurice) the SEF largely improves the $Q_{95}(|e|)$ of BOA. For example, for a 48-hour lead time in Chamonix, the RMSE and $Q_{95}(|e|)$ of $\text{BOA}^{s}$ are 1.54°C and 3.06°C, respectively, compared with 1.68°C and 3.53°C for BOA without the SEF

Montpellier, Orange, Perpignan and Tarbes are the only stations with a negative median and the SEF does not tend to improve the $Q_{95}(|e|)$ for these stations. In these stations, either the GBRT is not able to predict when BOA makes large errors, and/or the specialized experts are not specialized enough, i.e. they make predictions which are too cold or too warm when they are woken up.

In order to estimate which of the two adds more noise, we made two oracle strategies. The first one was $\text{BOA}^{s}$ combined with a theoretical GBRT making perfect $\widehat{e}_t'$ predictions. The second one was $\text{BOA}^{s}$ which suffers no loss when the specialized experts are woken up by our GBRT because we considered the specialized experts to be perfectly adapted. These two oracles have an RMSE across all the stations and lead times of 1.08°C and 0.89°C respectively and a $Q_{95}(|e|)$ of 2.12°C and 1.95°C, suggesting that the major improvement can be made if the sleeping experts are more specialized.

\subsection{Special case studies}
\label{special_case_study}

\subsubsection{Behavior around $0$°C}
A difficult and critical weather situation to forecast is when during the winter - where lots of weather systems are moving through France - the temperature is close to $0$°C. In this case, it is difficult to predict if it will snow or rain. That is why we computed (over all stations and lead times) the RMSE and the $Q_{95}(|e|)$ only with observations from November to December ranged between -1°C and +1°C (20852 observations). For these observations, the RMSE and $Q_{95}(|e|)$ are 1.29°C and 2.15°C respectively for $\text{BOA}^s$ which is better than 1.45°C and 2.53°C for BOA. Hence, in this kind of high-stakes situations, $\text{BOA}^s$ outperformed BOA.

\subsubsection{Several consecutive large errors}
\label{consecutive_errors}
\begin{figure*} [h]
 \centerline{\includegraphics[width=39pc]{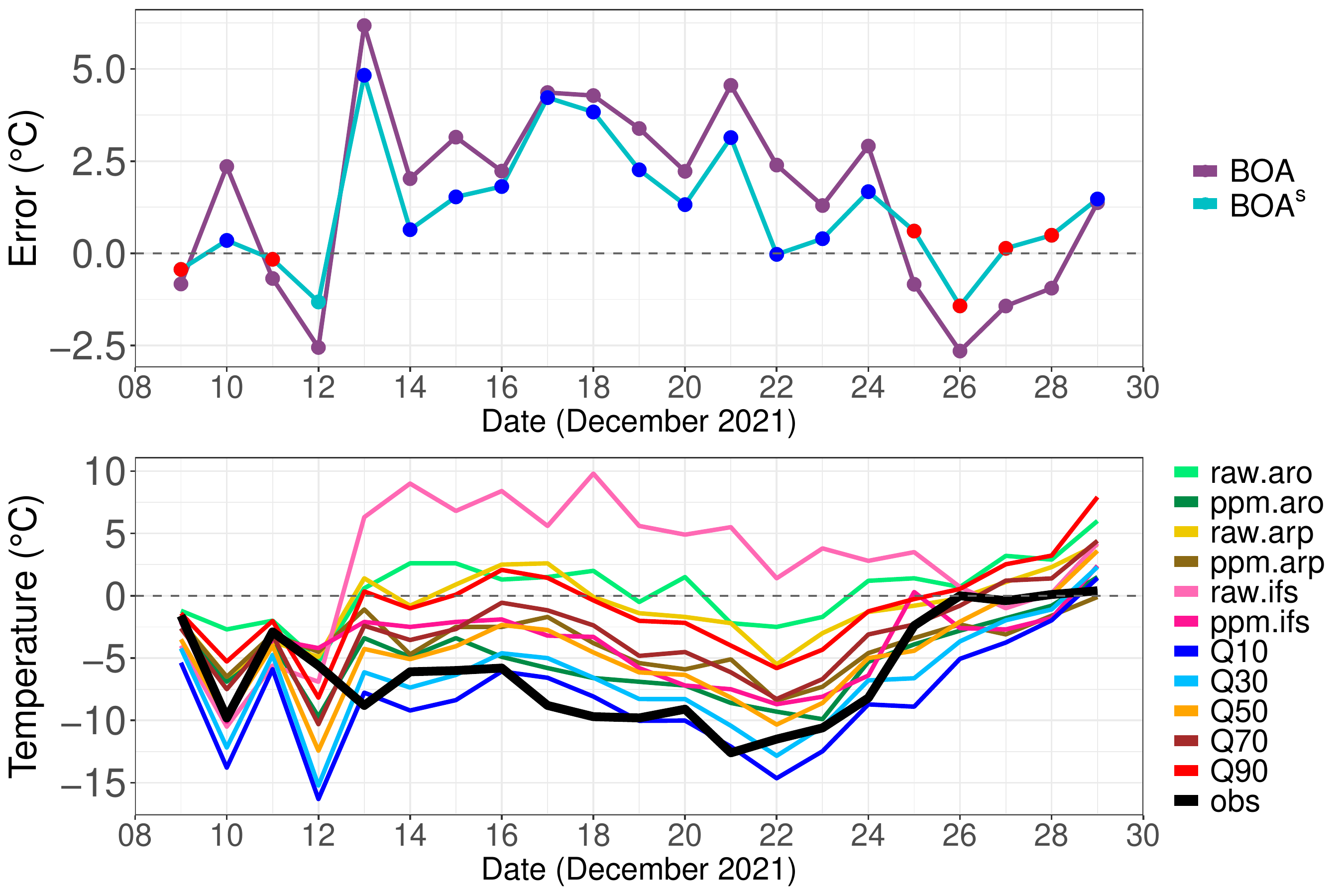}}
  \caption{Event of December 2021 at Chamonix (48-hour lead time). (Top) Error of Bernstein Online Aggregation (BOA) (purple) and  BOA in the Sleeping expert framework ($\text{BOA}^s$) (turquoise): blue points on the turquoise curve indicate that the Q10 and Q30 experts are awake (Q70 and Q90 asleep), red points indicate the opposite, and turquoise points indicate that all biased experts are asleep. The closer the curves are to the $y=0$ axis (dashed gray line), the better the performance. (Bottom) Experts predictions in December 2021 and in black the observed temperature. The experts are the "raw" and Post-Processed Models (ppm) AROME (aro), IFS (ifs), ARPEGE (arp) and the post-processed quantiles of the PEARP model (Q10, Q30, Q50, Q70, Q90).}\label{BOAvsBOA_sleep}
\end{figure*}

Now let us look at the Chamonix station, where temperatures are difficult to predict due to the mountainous terrain, where katabatic winds and radiative cooling effects remain poorly captured by the raw NWP models and PPMs. In Figure \ref{BOAvsBOA_sleep}, we can see the temperature observations of December 2021 along with the 48-hour lead time predictions of the different experts and the prediction errors of $\text{BOA}^{s}$ and BOA in Chamonix.

We consider this event as extreme because, as shown in the lower plot of Figure \ref{BOAvsBOA_sleep}, all experts (except the specialized Q10 and Q30 experts) produce forecasts that are far too warm throughout the entire period.

It appears that BOA, without the SEF, does not react to this event, since BOA is not able to quickly put more weight on the Q10 and Q30 experts as one can see in Figure \ref{BOA_sleep_BOA_window1253_ech48_74056001_2020-03-30_2023-09-03_weights} around iteration 630 (between the dashed grey lines) which corresponds to this event. And BOA follows more or less the predictions of the best expert so far which is the ppm.aro (ppm.aro has the largest weight amongst the experts in Figure \ref{BOA_sleep_BOA_window1253_ech48_74056001_2020-03-30_2023-09-03_weights}). So for BOA it is just "business as usual". Indeed, without the SEF, the event is far too short relative to the number of past iterations over which the experts have already accumulated losses, and the Q10 and Q30 have no time to recover large weights.

On the upper plot of Figure \ref{BOAvsBOA_sleep}, the points on the turquoise line representing $\text{BOA}^{s}$'s error, are blue when the Q10 and Q30 experts are awake and the Q70 and Q90 experts asleep (vice versa for the red points) and turquoise when all biased experts are asleep. One can notice that for this event, GBRT always wakes up the suitable biased experts (e.g. Q10 and Q30 when the error of BOA is positive) in order to bring $\text{BOA}^{s}$'s error closer to zero. This explains why $\text{BOA}^{s}$'s turquoise error line is closer to the dashed grey line and why this event is much better predicted by $\text{BOA}^{s}$ than by $\text{BOA}$.

Indeed, the error of $\text{BOA}^{s}$ is almost always closer to the dashed grey line hence smaller than BOA's error. This is due to the activation of the sleeping experts Q10, Q30, Q70 and Q90 at the right time.

\section{Conclusion}
We believe that the major improvement due to the SEF lies in the fact that it made the BOA strategy more reactive without adding noise, on the contrary it outperformed BOA with regard to the RMSE and the $Q_{95}(|e|)$. $\text{BOA}^{s}$'s RMSE and $Q_{95}(|e|)$ are respectively statistically significantly lower than BOA's RMSE and $Q_{95}(|e|)$ for 36.4\% and 5.9\% of the station - lead time pairs.

$\text{BOA}^{s}$ even significantly outperformed the best fixed convex combination of experts in hindsight for the RMSE and the $Q_{95}(|e|)$, something that BOA had not been able to achieve.

We observed that the contribution of the SEF to EA depends on the station and the lead time, sometimes largely improving the forecasts and sometimes adding some noise.

The SEF can help EA to largely improve forecasts of stations where temperature prediction is difficult such as in mountains. In particular the predictions of the December event in Chamonix - presented in the special case study - for which BOA made several consecutive large prediction errors were improved by the SEF.  The SEF also seems to be more accurate in critical situations where it is difficult to predict whether it will rain or snow because the temperature is close to $0$°C.

One could argue that retraining the GBRT at every iteration is a drawback of our method. Yet, we believe that this is a strength of our method. Indeed, the method is online and adaptive and does not need much maintenance if it is operational, even when the raw NWP models change.

Furthermore, the training task of the GBRT is fast enough to allow the use of our method operationally since we need less than half an hour to make the 1253 predictions for one station and one lead time on a supercomputer node.

Additionally, we would like, as in Part I, to highlight the robustness of EA and our method, since one could replace any of our experts by any temperature predictions of its choice, or add any temperature forecast to our set of experts especially the temperature predictions of the fully machine learning models such as in \citet{bi_pangu-weather_2022,nguyen_climax_2023,lang_aifs_2024}.

Concerning further studies, there are two main ways to improve the SEF. The first way is to improve the specialized experts when they are awake, or to add other sleeping experts that were not available for this study, such as the 5\% and 95\% quantiles of an ensemble NWP model.  The second way is to improve the GBRT. For example, one could try to use other variables such as wind direction or sea level pressure to improve the GBRT. It would also be interesting to replace the GBRT by other methods, such as neural networks or simpler machine learning algorithms with fewer (hyper)parameters. Indeed, one drawback of GBRT is the large number of hyperparameters to be tuned. Finally, the framework of experts that report their confidence \citep{blum_external_2007} could present an opportunity to further improve the EA by awakening the experts by an amount correlated to their quality for the current event.

We believe that the SEF is more adapted to temperature predictions than strategies designed to compete with sequences of experts as in \citet{herbster_tracking_1998,mourtada_efficient_2017}, because of the nature of the changes we want to catch and which the standard EA are unable to catch. Indeed, these changes are rapid changes of the best expert, very often events where a previously bad expert becomes good for one or for few iterations. The FS EA of \citet{herbster_tracking_1998} seems more adapted to deal with abrupt changes that last quite a long time, which seems unlikely in our case. Indeed, it is unrealistic that at some point one really bad expert becomes, for a long period of time, really good and vice versa. FS is more reactive and can put more weight on the biased experts than the EAs of Part I but our experiments with FS added noise compared to all the EAs of Part I.

\clearpage
\section*{Acknowledgments}
The authors would like to thank the three anonymous reviewers for their helpful comments and suggestions that greatly improved the quality of the article. The authors would like to acknowledge the support of the French Agence Nationale de la Recherche (ANR) under reference ANR20-CE40-0025-01 (T-REX project). Finally, the authors are grateful to Hermann Pfitzner for his thorough review.

\section*{Data availability}
The data and the code used for the experiments are available at https://github.com/pfitznerl/agregation\_2025.

\section*{Appendix: List of the variables used in the GBRT}

\begin{itemize}
    \item The observed temperature at the first lead time of the NWP models and the differences between the first lead time observation and the corresponding prediction of the aggregation and the experts.
    \item Standard deviation of the set of experts and standard deviation of the PEARP quantiles.
    \item Difference between the mean of the set of experts and the minimum and maximum of the set of experts.
    \item The difference between the mean of the experts and the prediction of the aggregation. Under certain conditions (the ensemble has to be calibrated), the error of the mean of an ensemble is correlated to the variance of the ensemble \citep{fortin_why_2014}.
    \item Across one run, the mean of the variance of the experts. Across one run, the variance of the variance of the experts. And the difference between the variance at the lead time that we want to predict and the mean variance over the lead times.
    \item The predictions and the standard deviation predictions of a Viking kalman filter \citep{de_vilmarest_viking_2021,de_vilmarest_state-space_2022} applied to AROME if possible or ARPEGE otherwise. We did not add the kalman filters applied to the other experts because the correlation between the kalman filters was considerably high.
\end{itemize}

\section*{Appendix: Confusion matrix and Equitable Skill Score}
\label{gerrity_skill_score}

If there are $K\in \mathbb{N}$ possible classes/observations, and one has the contingency table with $n_{\widehat{e},e}\in \mathbb{R}$ on row $\widehat{e}=1,\ldots,K$ and column $e=1,\ldots,K$ (see Table \ref{contingency_table} for the contingency table corresponding to our study with 3 classes), then the corresponding ESS (often called Gerrity skill score) is:
\begin{equation}
    ESS=\sum_{\widehat{e}=1}^K\sum_{e=1}^K p_{\widehat{e},e}s_{\widehat{e},e}
\end{equation}

where $p_{\widehat{e},e}=n_{\widehat{e},e} / (\sum_{r=1}^K\sum_{s=1}^K n_{r,s})$ is the joint distribution of prediction $\widehat{e}$ and observation $e$ and $s_{\widehat{e},e}$ is the reward for this event, of the scoring matrix with coefficient $s_{\widehat{e},e}$ on line $\widehat{e}$ and column $e$, $1 \leqslant \widehat{e},e \leqslant K$.
\citet{gerrity_note_1992} provides a way in which to compute the scoring matrix in order to have an equitable skill score:
\begin{align}
    s_{\widehat{e},e}=\frac{1}{K-1} \left( \sum_{r=1}^{\widehat{e}-1} a_r^{-1} - \sum_{r=\widehat{e}}^{e-1} 1 + \sum_{r=e}^{K-1} a_r \right),
\end{align}

with $a_{e}=\frac{1-\sum_{r=1}^e p_{e}}{\sum_{r=1}^e p_{e}}$ and $p_e=\sum_{\widehat{e}=1}^K p_{\widehat{e},e}$.

\citet{gerrity_note_1992} also showed that the mean of the $K-1$ Peirce Skill Score (PSS) associated to the two-class problems generated by partitioning the original contingency table at its $K-1$ thresholds is also an equitable skill score:
\begin{align}
    ESS=\frac{1}{K-1} \sum_{r=1}^{K-1} PSS(r),
\end{align}

where for $r=1,\ldots,K$, PSS(r) is the PSS associated to the $r^{th}$ threshold with:
\begin{align}
    PSS(r)=\frac{ \left( \sum_{e=1}^r \sum_{\widehat{e}=1}^r p_{\widehat{e},e} \right) \left( \sum_{e=r+1}^K \sum_{\widehat{e}=r+1}^K p_{\widehat{e},e} \right) - \left( \sum_{e=1}^r \sum_{\widehat{e}=r+1}^K p_{\widehat{e},e} \right) \left( \sum_{e=r+1}^K \sum_{\widehat{e}=1}^r p_{\widehat{e},e} \right) }{ \left( \left( \sum_{e=1}^r \sum_{\widehat{e}=1}^r p_{\widehat{e},e} \right) + \left( \sum_{e=r+1}^K \sum_{\widehat{e}=1}^r p_{\widehat{e},e} \right) \right) \left( \left( \sum_{e=1}^r \sum_{\widehat{e}=r+1}^K p_{\widehat{e},e} \right) + \left( \sum_{e=r+1}^K \sum_{\widehat{e}=r+1}^K p_{\widehat{e},e} \right) \right) }.
\end{align}

An ESS of 0 indicates that the model is like a random prediction. A strictly negative (positive) ESS shows that the model is worse (better) than a random prediction according to the climatology, and an ESS of 1 (-1) indicates that the model is always right (wrong).

\begin{table*}[h]
\begin{center}
\begin{tabular}{c|ccc}
 & \multicolumn{3}{c}{\textbf{Observations}}\\
\hline
\textbf{Predictions} & $e_t \leqslant -e^\text{thr}$ & $-e^\text{thr} < e_t < e^\text{thr}$ & $e_t \geqslant e^\text{thr}$  \\
\hline
$e_t$ $\leqslant -e^\text{thr}$ & $n_{1,1}$ & $n_{1,2}$ & $n_{1,3}$ \\
$-e^\text{thr} < e_t < e^\text{thr}$ & $n_{2,1}$ & $n_{2,2}$ & $n_{2,3}$\\
$e_t \geqslant e^\text{thr}$ & $n_{3,1}$ & $n_{3,2}$ & $n_{3,3}$\\
\hline
\end{tabular}
\end{center}
\caption{Confusion matrix.}\label{contingency_table}
\end{table*}

\bibliographystyle{unsrtnat}
\bibliography{references.bib}  

@STRING{AN        = "Astrophys.\ Norv."}

@STRING{MA        = "Meteor.\ Appl."}

@inproceedings{gaillard_second-order_2014,
	title = {A {Second}-order {Bound} with {Excess} {Losses}},
	booktitle = {Conference on {Learning} {Theory}.},
	publisher = {PMLR},
	author = {Gaillard, Pierre and Stoltz, Gilles and van Erven, Tim},
	year = {2014},
	pages = {176--196},
}

@inproceedings{jun_improved_2017,
	title = {Improved {Strongly} {Adaptive} {Online} {Learning} using {Coin} {Betting}},
	booktitle = {Proceedings of the 20th {International} {Conference} on {Artificial} {Intelligence} and {Statistics}},
	publisher = {PMLR},
	author = {Jun, Kwang-Sung and Orabona, Francesco and Wright, Stephen and Willett, Rebecca},
	year = {2017},
	pages = {943--951},
}

@article{vyugin_online_2019,
	title = {Online aggregation of unbounded losses using shifting experts with confidence},
	volume = {108},
	doi = {10.1007/s10994-018-5751-z},
	journal = {Machine Learning},
	author = {V’yugin, Vladimir and Trunov, Vladimir},
	year = {2019},
	pages = {425--444},
}

@article{taillardat_calibrated_2016,
	title = {Calibrated {Ensemble} {Forecasts} {Using} {Quantile} {Regression} {Forests} and {Ensemble} {Model} {Output} {Statistics}},
	volume = {144},
	issn = {1520-0493, 0027-0644},
	url = {https://journals.ametsoc.org/view/journals/mwre/144/6/mwr-d-15-0260.1.xml},
	doi = {10.1175/MWR-D-15-0260.1},
	journal = {Monthly Weather Review},
	author = {Taillardat, Maxime and Mestre, Olivier and Zamo, Michaël and Naveau, Philippe},
	year = {2016},
	pages = {2375--2393},
}

@inproceedings{mourtada_efficient_2017,
	title = {Efficient tracking of a growing number of experts},
	booktitle = {Proceedings of the 28th {International} {Conference} on {Algorithmic} {Learning} {Theory}},
	publisher = {PMLR},
	author = {Mourtada, Jaouad and Maillard, Odalric-Ambrym},
    pages={517-539},
	year = {2017},
}

@article{herbster_tracking_1998,
	title = {Tracking the {Best} {Expert}},
	volume = {32},
	issn = {1573-0565},
	url = {https://doi.org/10.1023/A:1007424614876},
	doi = {10.1023/A:1007424614876},
	language = {en},
	urldate = {2021-01-11},
	journal = {Machine Learning},
	author = {Herbster, Mark and Warmuth, Manfred K.},
	month = aug,
	year = {1998},
	pages = {151--178},
}

@article{gaillard_opera_2016,
    title={Opera: Online prediction by expert aggregation.},
    journal={ R package version 1.1.1},
    author={Gaillard, Pierre and Goude, Yannig},
    year={2016},
    url = {https://cran.r-project.org/src/contrib/Archive/opera/opera_1.1.1.tar.gz}
}

@inproceedings{freund_using_1997,
	title = {Using and combining predictors that specialize},
	doi = {10.1145/258533.258616},
	booktitle = {Proceedings of the twenty-ninth annual {ACM} symposium on {Theory} of computing},
	author = {Freund, Yoav and Schapire, Robert E. and Singer, Yoram and Warmuth, Manfred K.},
	adress = {New York, NY, USA},
    organization={Association for Computing Machinery},
	year = {1997},
	pages = {334--343}
}

@article{devaine_forecasting_2013,
	title = {Forecasting electricity consumption by aggregating specialized experts},
	volume = {90},
	issn = {1573-0565},
	url = {https://doi.org/10.1007/s10994-012-5314-7},
	doi = {10.1007/s10994-012-5314-7},
	language = {en},
	urldate = {2021-01-11},
	journal = {Machine Learning},
	author = {Devaine, Marie and Gaillard, Pierre and Goude, Yannig and Stoltz, Gilles},
	month = feb,
	year = {2013},
	pages = {231--260},
}

@book{cesa-bianchi_prediction_2006,
	title = {Prediction, {Learning}, and {Games}},
	publisher = {Cambridge University Press},
	author = {Cesa-Bianchi, Nicolo and Lugosi, Gabor},
	month = mar,
	year = {2006},
    pages={407}
}

@article{blum_external_2007,
	title = {From {External} to {Internal} {Regret}},
	volume = {8},
	journal = {Journal of Machine Learning Research},
	author = {Blum, Avrim and Mansour, Yishay},
	year = {2007},
	pages = {1307--1324},
}

@article{blum_empirical_1997,
	title = {Empirical {Support} for {Winnow} and {Weighted}-{Majority} {Algorithms}: {Results} on a {Calendar} {Scheduling} {Domain}},
	volume = {26},
	doi = {10.1023/A:1007335615132},
	urldate = {2021-01-11},
	journal = {Machine Learning},
	author = {Blum, Avrim},
	year = {1997},
	pages = {5--23},
}

@article{silva_using_2022,
	title = {Using an {Explainable} {Machine} {Learning} {Approach} to {Characterize} {Earth} {System} {Model} {Errors}: {Application} of {SHAP} {Analysis} to {Modeling} {Lightning} {Flash} {Occurrence}},
	volume = {14},
	doi = {10.1029/2021MS002881},
	journal = {Journal of Advances in Modeling Earth Systems},
	author = {Silva, Sam J. and Keller, Christoph A. and Hardin, Joseph},
	year = {2022},
	pages = {e2021MS002881},
}

@article{taillardat_research_2020,
	title = {From research to applications – examples of operational ensemble post-processing in {France} using machine learning},
	volume = {27},
	doi = {10.5194/npg-27-329-2020},
	journal = {Nonlinear Processes in Geophysics},
	author = {Taillardat, Maxime and Mestre, Olivier},
	month = may,
	year = {2020},
	pages = {329--347},
}

@article{rasp_neural_2018,
	title = {Neural {Networks} for {Postprocessing} {Ensemble} {Weather} {Forecasts}},
	volume = {146},
	doi = {10.1175/MWR-D-18-0187.1},
	journal = {Monthly Weather Review},
	author = {Rasp, Stephan and Lerch, Sebastian},
	month = nov,
	year = {2018},
	pages = {3885--3900},
}

@misc{bi_pangu-weather_2022,
	title = {Pangu-{Weather}: {A} {3D} {High}-{Resolution} {Model} for {Fast} and {Accurate} {Global} {Weather} {Forecast}},
	doi = {10.48550/arXiv.2211.02556},
	publisher = {arXiv},
	author = {Bi, Kaifeng and Xie, Lingxi and Zhang, Hengheng and Chen, Xin and Gu, Xiaotao and Tian, Qi},
	year = {2022},
}

@misc{nguyen_climax_2023,
	title = {{ClimaX}: {A} foundation model for weather and climate},
	language = {en},
	publisher = {arXiv},
	author = {Nguyen, Tung and Brandstetter, Johannes and Kapoor, Ashish and Gupta, Jayesh K. and Grover, Aditya},
	year = {2023},
    doi={10.48550/arXiv.2301.10343},
}

@misc{de_vilmarest_viking_2021,
	title = {Viking: {Variational} {Bayesian} {Variance} {Tracking}},
	shorttitle = {Viking},
	doi = {10.48550/arXiv.2104.10777},
	urldate = {2024-05-28},
	publisher = {arXiv},
	author = {De Vilmarest, Joseph and Wintenberger, Olivier},
	year = {2021},
}

@article{de_vilmarest_state-space_2022,
	title = {State-{Space} {Models} for {Online} {Post}-{Covid} {Electricity} {Load} {Forecasting} {Competition}},
	volume = {9},
	doi = {10.1109/OAJPE.2022.3141883},
	journal = {IEEE Open Access Journal of Power and Energy},
	author = {De Vilmarest, Joseph and Goude, Yannig},
	year = {2022},
	pages = {192--201},
}

@inproceedings{adamskiy_closer_2012,
	address = {Berlin, Heidelberg},
	series = {Lecture {Notes} in {Computer} {Science}},
	title = {A {Closer} {Look} at {Adaptive} {Regret}},
	doi = {10.1007/978-3-642-34106-9_24},
	booktitle = {Algorithmic {Learning} {Theory}},
	publisher = {Springer},
	author = {Adamskiy, Dmitry and Koolen, Wouter M. and Chernov, Alexey and Vovk, Vladimir},
	year = {2012},
	pages = {290--304},
}

@article{wintenberger_stochastic_2024,
	title = {Stochastic online convex optimization. {Application} to probabilistic time series forecasting},
	volume = {18},
	doi = {10.1214/23-EJS2208},
	journal = {Electronic Journal of Statistics},
	author = {Wintenberger, Olivier},
	year = {2024},
	pages = {429--464},
}

@misc{lang_aifs_2024,
	title = {{AIFS} - {ECMWF}'s data-driven forecasting system},
	doi = {10.48550/arXiv.2406.01465},
	publisher = {arXiv},
	author = {Lang, Simon and Alexe, Mihai and Chantry, Matthew and Dramsch, Jesper and Pinault, Florian and Raoult, Baudouin and Clare, Mariana C. A. and Lessig, Christian and Maier-Gerber, Michael and Magnusson, Linus and Bouallègue, Zied Ben and Nemesio, Ana Prieto and Dueben, Peter D. and Brown, Andrew and Pappenberger, Florian and Rabier, Florence},
	year = {2024},
}

@inproceedings{chen_xgboost_2016,
	title = {{XGBoost}: {A} {Scalable} {Tree} {Boosting} {System}},
	doi = {10.1145/2939672.2939785},
	booktitle = {Proceedings of the 22nd {ACM} {SIGKDD} {International} {Conference} on {Knowledge} {Discovery} and {Data} {Mining}},
	publisher = {Association for Computing Machinery},
	author = {Chen, Tianqi and Guestrin, Carlos},
	year = {2016},
	pages = {785--794},
}

@inproceedings{hazan_adaptive_2007,
  title={Adaptive algorithms for online decision problems},
  author={Hazan, Elad and Seshadhri, Comandur},
  booktitle={Electronic colloquium on computational complexity (ECCC)},
  year={2007}
}

@article{fei_factors_2020,
	title = {Factors {Affecting} the {Weakening} {Rate} of {Tropical} {Cyclones} over the {Western} {North} {Pacific}},
	volume = {148},
	doi = {10.1175/MWR-D-19-0356.1},
	journal = {Monthly Weather Review},
	author = {Fei, Rong and Xu, Jing and Wang, Yuqing and Yang, Chi},
	year = {2020},
	pages = {3693--3712}
}

@article{potdar_toward_2021,
	title = {Toward {Predicting} {Flood} {Event} {Peak} {Discharge} in {Ungauged} {Basins} by {Learning} {Universal} {Hydrological} {Behaviors} with {Machine} {Learning}},
	volume = {22},
	doi = {10.1175/JHM-D-20-0302.1},
	journal = {Journal of Hydrometeorology},
	author = {Potdar, Akhil Sanjay and Kirstetter, Pierre-Emmanuel and Woods, Devon and Saharia, Manabendra},
	year = {2021},
	pages = {2971--2982},
}

@article{mony_evaluating_2021,
	title = {Evaluating {Foehn} {Occurrence} in a {Changing} {Climate} {Based} on {Reanalysis} and {Climate} {Model} {Data} {Using} {Machine} {Learning}},
	volume = {36},
	doi = {10.1175/WAF-D-21-0036.1},
	journal = {Weather and Forecasting},
	author = {Mony, Christoph and Jansing, Lukas and Sprenger, Michael},
	year = {2021},
	pages = {2039--2055},
}

@article{flora_using_2021,
	title = {Using {Machine} {Learning} to {Generate} {Storm}-{Scale} {Probabilistic} {Guidance} of {Severe} {Weather} {Hazards} in the {Warn}-on-{Forecast} {System}},
	volume = {149},
	doi = {10.1175/MWR-D-20-0194.1},
	journal = {Monthly Weather Review},
	author = {Flora, Montgomery L. and Potvin, Corey K. and Skinner, Patrick S. and Handler, Shawn and McGovern, Amy},
	month = may,
	year = {2021},
	pages = {1535--1557},
}

@article{liu_correction_2022,
	title = {Correction of {Overestimation} in {Observed} {Land} {Surface} {Temperatures} {Based} on {Machine} {Learning} {Models}},
	volume = {35},
	doi = {10.1175/JCLI-D-21-0447.1},
	journal = {Journal of Climate},
	author = {Liu, Fa and Wang, Xunming and Sun, Fubao and Wang, Hong and Wu, Lifeng and Zhang, Xuanze and Liu, Wenbin and Che, Huizheng},
	year = {2022},
	pages = {5359--5377},
}

@article{nanda_soil_2020,
	title = {Soil {Temperature} {Dynamics} at {Hillslope} {Scale}—{Field} {Observation} and {Machine} {Learning}-{Based} {Approach}},
	volume = {12},
	doi = {10.3390/w12030713},
	journal = {Water},
	author = {Nanda, Aliva and Sen, Sumit and Sharma, Awshesh Nath and Sudheer, K. P.},
	year = {2020},
	pages = {713},
}

@article{ma_prediction_2020,
	title = {Prediction of outdoor air temperature and humidity using {Xgboost}},
	volume = {427},
	doi = {10.1088/1755-1315/427/1/012013},
	journal = {IOP Conference Series: Earth and Environmental Science},
	author = {Ma, Xiaoming and Fang, Cong and Ji, Junping},
	year = {2020},
	pages = {012013},
}

@article{gerrity_note_1992,
	title = {A {Note} on {Gandin} and {Murphy}'s {Equitable} {Skill} {Score}},
	volume = {120},
	doi = {10.1175/1520-0493(1992)120<2709:ANOGAM>2.0.CO;2},
	journal = {Monthly Weather Review},
	author = {Gerrity, Joseph P.},
	year = {1992},
	pages = {2709--2712},
}

@article{wmo_manual_2010,
	title = {Manual on the {Global} {Data}-{Processing} and {Forecasting} {System}. {Attachment} {II}, {Global} Aspects},
	journal = {World Meteorological Organization},
	author = {WMO},
	year = {2010}
}

@article{fortin_why_2014,
	title = {Why {Should} {Ensemble} {Spread} {Match} the {RMSE} of the {Ensemble} {Mean}?},
	volume = {15},
	doi = {10.1175/JHM-D-14-0008.1},
	urldate = {2023-09-01},
	journal = {Journal of Hydrometeorology},
	author = {Fortin, V. and Abaza, M. and Anctil, F. and Turcotte, R.},
	year = {2014},
	pages = {1708--1713}
}

@article{zhao_linear_2024,
	title = {Linear discriminant analysis},
	volume = {4},
	issn = {2662-8449},
	url = {https://www.nature.com/articles/s43586-024-00346-y},
	doi = {10.1038/s43586-024-00346-y},
	journal = {Nature Reviews Methods Primers},
	author = {Zhao, Shuping and Zhang, Bob and Yang, Jian and Zhou, Jianhang and Xu, Yong},
	year = {2024},
	note = {Publisher: Nature Publishing Group},
	pages = {70},
}

@article{marzban_neural_2003,
	title = {Neural {Networks} for {Postprocessing} {Model} {Output}: {ARPS}},
	volume = {131},
	url = {https://journals.ametsoc.org/view/journals/mwre/131/6/1520-0493_2003_131_1103_nnfpmo_2.0.co_2.xml},
	doi = {10.1175/1520-0493(2003)131<1103:NNFPMO>2.0.CO;2},
	journal = {Monthly Weather Review},
	author = {Marzban, Caren},
	year = {2003},
	pages = {1103--1111},
}

@article{diebold_comparing_1995,
	title = {Comparing {Predictive} {Accuracy}},
	volume = {20},
	doi = {10.1198/073500102753410444},
	journal = {Journal of Business \& Economic Statistics},
	author = {Diebold, Francis X and Mariano, Robert S},
	year = {1995},
	pages = {134--144},
}

@article{benjamini_controlling_1995,
	title = {Controlling the {False} {Discovery} {Rate}: {A} {Practical} and {Powerful} {Approach} to {Multiple} {Testing}},
	volume = {57},
	doi = {10.1111/j.2517-6161.1995.tb02031.x},
	journal = {Journal of the Royal Statistical Society: Series B (Methodological)},
	author = {Benjamini, Yoav and Hochberg, Yosef},
	year = {1995},
	pages = {289--300},
}

@article{johnson_two-sample_1987,
	title = {Two-{Sample} {Rank} {Tests} for {Detecting} {Changes} {That} {Occur} in a {Small} {Proportion} of the {Treated} {Population}},
	volume = {43},
	doi = {10.2307/2532001},
	journal = {Biometrics},
	author = {Johnson, Richard A. and Verrill, Steve and Moore, Dan H.},
	year = {1987},
	pages = {641--655},
}

@article{monteleoni_tracking_2011,
	title = {Tracking climate models},
	volume = {4},
	doi = {10.1002/sam.10126},
	journal = {Statistical Analysis and Data Mining: The ASA Data Science Journal},
	author = {Monteleoni, Claire and Schmidt, Gavin A. and Saroha, Shailesh and Asplund, Eva},
	year = {2011},
	pages = {372--392},
}






\end{document}